\documentclass[review,12pt]{elsarticle}

\usepackage{bm}  
\usepackage{amssymb}
\usepackage{amsthm}
\usepackage{amsmath}
\usepackage[linesnumbered, ruled, lined]{algorithm2e}
\graphicspath{{}} 
\usepackage{booktabs}
\usepackage{multirow}
\usepackage{verbatim}
\usepackage{hyperref}
\usepackage{lipsum}
\usepackage{graphicx}
\usepackage{subfigure}

\usepackage[outline]{contour}
\newcommand*{\fancy}[1]{{\color{white}\contour{black}{#1}}}

\usepackage{nomencl}
\makenomenclature

\renewcommand{\nompreamble}{The next list describes the symbols that are used within the body of the document.}

\journal{Elsevier Journal}

\begin{document}
\begin{frontmatter}



\title{An SPH formulation for general plate and shell structures with finite deformation and large rotation}
\author[myfirstaddress]{Dong Wu}
\ead{dong.wu@tum.de}
\author[myfirstaddress,mysecondaddress]{Chi Zhang}
\ead{c.zhang@tum.de}
\author[myfirstaddress]{Xiangyu Hu\corref{mycorrespondingauthor}}
\cortext[mycorrespondingauthor]{Corresponding author.}
\ead{xiangyu.hu@tum.de}
\address[myfirstaddress]{TUM School of Engineering and Design, Technical University of Munich, 85748 Garching, Germany}
\address[mysecondaddress]{Huawei Technologies Munich Research Center, 80992 Munich, Germany}

\begin{abstract}
In this paper, 
we propose a reduced-dimensional smoothed particle hydrodynamics (SPH) formulation 
for quasi-static and dynamic analyses of plate and shell structures 
undergoing finite deformation and large rotation. 
By exploiting Uflyand–Mindlin plate theory, 
the present surface-particle formulation is able to resolve
the thin structures by using only one layer of particles at the mid-surface.
To resolve the geometric non-linearity 
and capture finite deformation and large rotation,
two reduced-dimensional linear-reproducing correction matrices are introduced, 
and weighted non-singularity conversions between the rotation angle and pseudo normal 
are formulated. 
A new non-isotropic Kelvin-Voigt damping is proposed especially 
for the both thin and moderately thick plate and shell structures 
to increase the numerical stability. 
In addition, 
a shear-scaled momentum-conserving hourglass control algorithm with 
an adaptive limiter 
is introduced to suppress the mismatches
between the particle position and pseudo normal 
and those estimated with the deformation gradient. 
A comprehensive set of test problems,
for which the analytical or numerical results from literature  
or those of the volume-particle SPH model
are available for quantitative and qualitative comparison,
are examined to demonstrate the accuracy and stability of the present method.
\end{abstract}



\end{frontmatter}


\section{Introduction}
For computational continuum dynamics, 
as alternatives to conventional mesh-based methods, 
e.g. finite element method (FEM) and finite volume method (FVM), 
meshless methods have flourished in the past decades 
\cite{belytschko1996meshless, liu2005introduction, liu2010smoothed, zhang2022smoothed}. 
Smoothed particle hydrodynamics (SPH), 
initially developed by Lucy \cite{lucy1977numerical} 
and Gingold and Monaghan \cite{gingold1977smoothed} 
for astrophysical simulations, is one typical example. 
In SPH, 
the continuum is modeled by particles
associated with physical properties such as mass and velocity,
and the governing equations are discretized 
in the form of particle interactions using a Gaussian-like kernel function
\cite{monaghan2005smoothed, liu2010smoothed, monaghan2012smoothed}.
Since a significant number of physical system abstractions
can be realized through particle interactions, 
SPH has been used to model multi-physical systems 
within a unified computational framework \cite{zhang2021sphinxsys}, 
which is able to achieve seamless monolithic, strong and conservative coupling
\cite{matthies2003partitioned, matthies2006algorithms}. 

To achieve such a unified computational framework, 
it is crucial to discretize all relevant physics equations 
using effective and efficient SPH methods. 
In the case of plate and shell structures
which are omnipresent thin structures in scientific and 
engineering fields such as shipbuilding \cite{cerik2019simulation, peng2019meshfree}, 
aerospace \cite{totaro2009optimal}, 
and medical treatment \cite{laubrie2020new}, etc.,
the traditional full-dimensional or volume-particle SPH method, 
is not computationally efficient \cite{cleary1998modelling}. 
Since there are well-developed and matured 
reduced-dimensional theories,
such as Kirchhoff-Love \cite{love1887small} 
and Uflyand-Mindlin (or Mindlin-Reissner)  \cite{uflyand1948wave, mindlin1951influence, elishakoff2017vibrations, elishakoff2020handbook},
for plate and shell structures 
based on mid-surface reconstruction,
it is expected to develop 
the computationally much more efficient 
reduced-dimensional or surface-particle SPH method
with a single-layer of particles only. 

The early meshless methods for plates and shells
were based on Petrov or element-free Galerkin formulation
\cite{krysl1996analysis, li2000numerical, 
rabczuk2007meshfree, li2008numerical}, 
or the reproducing kernel particle method 
\cite{donning1998meshless, chen2006constrained, peng2018thick}. 
As for SPH, 
Maurel and Combescure \cite{maurel2008sph} first developed 
a surface-particle SPH method 
for total Lagrangian quasi-static and dynamic analyses of 
moderately thick plates and shells 
based on the Uflyand-Mindlin theory and the assumption of small deformation. 
In their work, 
besides an artificial viscosity term 
to alleviate numerical instability issues,
a stress point method is applied to temper hourglass or zero-energy modes
which exhibit in the traditional SPH method using collocated particles 
for both deformation and stress.
While being effective on preventing zero-energy modes, 
using stress points may faces several issues, 
such as how to locate or generate these points for complex geometries, 
complicated numerical algorithms 
and the compensation of computational efficiency 
\cite{ganzenmuller2015hourglass, wu2023essentially}.  
Nevertheless, this method was later applied in 
large deformation analyses  
by Ming et al . \cite{ming2013robust} 
and dynamic damage-fracture analyses 
by Caleyron et al. \cite{caleyron2012dynamic}.
Lin et al. \cite{lin2014efficient} 
developed a similar method for quasi-static analyses, 
but applied an artificial viscosity term 
based on membrane and shearing decomposition.
Ming et al. \cite{ming2015smoothed} 
first considered finite deformation by taking all strain terms into account
with the help of Gauss-Legendre quadrature 
for more accurately capturing of non-linear stress. 
In all of the surface-particle SPH methods mentioned above, 
the rotation angles of mid-surface 
are directly obtained from the pseudo normal in 
governing equations under the assumption of small rotation. 

In this work, 
we propose a collocated surface-particle SPH formulation 
for total Lagrangian quasi-static and dynamic analyses 
of general plate or shell structures, 
which may be thin or have moderate thickness,
involving finite deformation or/and large rotation. 
First, 
to better resolve the geometric non-linearity 
induced by finite deformation and large rotation, 
two new reduced-dimensional correction matrices
for linearly reproducing position and normal direction are introduced, 
and a weighted conversion algorithm,  
which achieves non-singularity under large rotation, is proposed. 
Second,
a new non-isotropic Kelvin-Voigt damping base on Ref. \cite{zhang2022artificial} 
is proposed for achieve good numerical stability 
for both thin and moderately thick plate or shell structures.
Third, 
in order to address hourglass modes using collocated particles only
other than introducing extra stress points,
drawing the inspiration from 
Refs. \cite{kondo2010suppressing, ganzenmuller2015hourglass},
a shear-scaled momentum-conserving formulation with an adaptive limiter
is developed by mitigating the discrepancy
between the actual particle position and pseudo normal 
and those estimated by the deformation gradient. 
A set of numerical examples 
involving quasi-static and dynamic analyses for both thin and moderately thick plate or shell structures are given. 
The results are compared with analytical, numerical solutions 
in literature or/and those obtained by the volume-particle SPH method
to demonstrate the numerical accuracy and stability of the present method. 

The remainder of this manuscript is organized as follows. 
Section \ref{sec:governingeq} introduces the theoretical model of plates and shells, including the kinematics, constitutive relation, stress correction and conservation equations. 
The proposed surface-particle SPH formulation, 
including the reduced-dimensional linear-producing correction matrices, 
weighted conversion algorithm, 
non-isotropic damping and momentum-conserving hourglass control,
is described in Section \ref{sec:SPH_method}.
Numerical examples are presented and discussed in Section \ref{sec:examples} 
and then concluding remarks are given in Section \ref{sec:conclusion}. 
For a better comparison and future opening for in-depth studies, 
all the computational codes of this work are released 
in the open-source repository of SPHinXsys 
\cite{zhang2020sphinxsys, zhang2021sphinxsys} at 
\texttt{https://github.com/Xiangyu-Hu/SPHinXsys}.
\nomenclature{$\left( \bullet \right)^0$}
{indicating the parameter $\left( \bullet \right)$ is defined at the initial configuration and global coordinate system}
\nomenclature{$\left( \bullet \right)^{0, L}$}
{indicating the parameter $\left( \bullet \right)$ is defined at the initial configuration and initial local coordinate system}
\nomenclature{$\left( \bullet \right)^l$}
{indicating the parameter $\left( \bullet \right)$ is defined at the current configuration and current local coordinate system}
\nomenclature{$\bm{X} = \left( X, Y, Z \right)$}
{global coordinate system}
\nomenclature{$\bm{\xi} = \left( \xi ,\eta ,\zeta \right)$ }
{initial local coordinate system}
\nomenclature{$\bm{x} = \left( x, y, z \right)$}
{current local coordinate system}
\section{Theoretical models}\label{sec:governingeq}
We first introduce 
the theoretical mode of 3D plate, 
and then that of 3D shell 
in which material points may possess different initial normal directions
leading to different initial local coordinate systems.
After that, we briefly describe the 2D plate and shell models, 
which resolve the plane strain problem, 
as a simplification of the 3D counterparts. 
\subsection{3D plate model}
We consider the Uflyand–Mindlin plate theory \cite{uflyand1948wave, mindlin1951influence} 
to account for transverse shear stress 
which is significant for moderately thick plates.
The theory implies that the plate behavior can be represented 
by one layer of material points at its mid-surface, 
as shown in Figure \ref{figs:discretization_plate}.  
\subsubsection{Kinematics}
\begin{figure}[htb!]
	\centering
	\includegraphics[trim = 2mm 8mm 2mm 2mm, width = 0.9\textwidth]{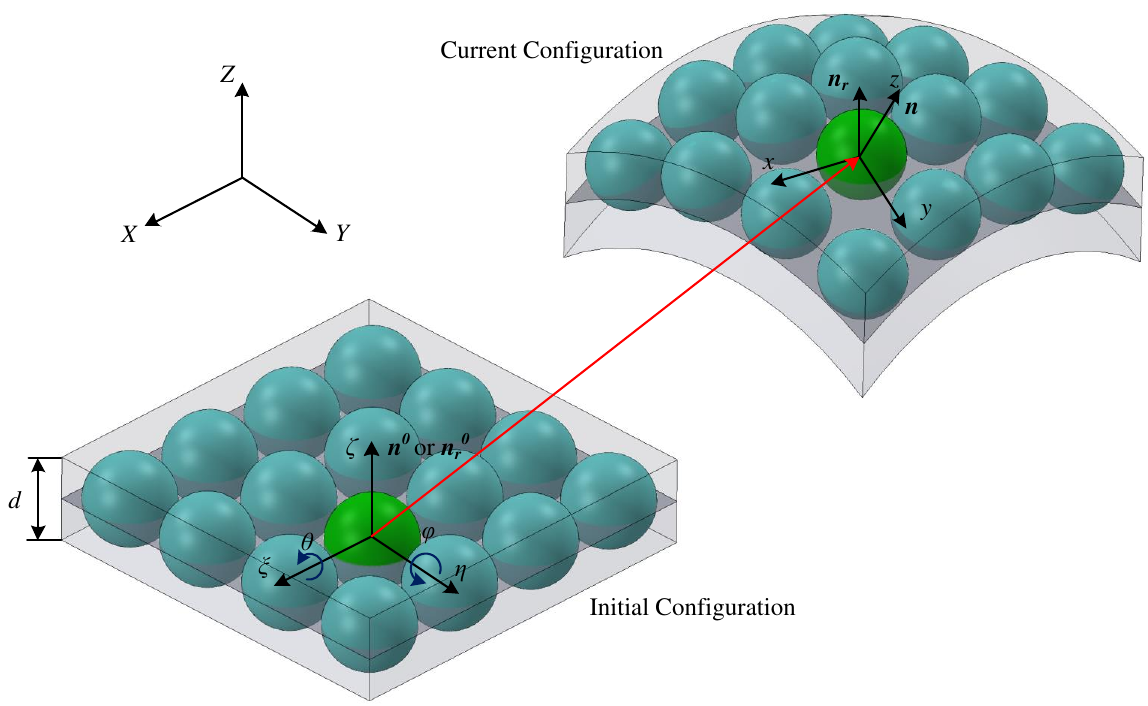}
	\caption{Schematic of a 3D plate model.}
	\label{figs:discretization_plate}
\end{figure}

We introduce $\bm{X} = \left( X, Y, Z \right)$ to represent the global coordinate system, 
and $\bm{\xi} = \left( \xi,\eta, \zeta \right)$ and $\bm{x} = \left( x, y, z \right)$, 
associated with so-called pseudo-normal vector $ \bm{n}$, 
to denote the initial and current local coordinate systems, respectively.
Note that the initial local coordinate system is same with the global one for plate.
Each material point possesses five degrees of freedom, viz.,
three  translations $ \bm{u}  = \left\{ u,v,w \right\}^{\operatorname{T}} $ 
and two rotations 
$\bm{\theta}  = \left\{ \theta, \varphi \right\}^{\operatorname{T}} $
expressed in the global coordinates.
Positive values of $\theta$ and $\varphi$ indicate 
that the plate is rotated anticlockwise around the coordinate axis
when the axis points toward the observer and the coordinate system is right-handed.
The two rotations are used to update the pseudo-normal
$ \bm{n} = \left\{ n_{1},n_{2},n_{3} \right\}^{\operatorname{T}}$ 
which is also defined in the global coordinate system 
and remains straight but is not necessarily perpendicular to the mid-surface, 
i.e., the pseudo normal may be different with the real normal $ \bm{n}_r$, 
as shown in Figure \ref{figs:discretization_plate}. 
Note that $ \bm{n}^0 = \left\{ 0, 0, 1 \right\}^{\operatorname{T}} $ 
denotes the pseudo-normal in the initial configuration
with the superscript $\left( \bullet \right)^0$ denoting 
the initial configuration.

For a 3D plate, 
the position $ \bm{r} $ of a material point at a distance $ \chi $ 
away from the mid-surface along the pseudo normal $\bm n$
can be expressed as 
\nomenclature{$ \bm{u}  = \left\{ u,v,w \right\}^{\operatorname{T}} $}{displacements}
\nomenclature{$\bm{\theta}  = \left\{ \theta, \varphi \right\}^{\operatorname{T}} $}{rotations}
\nomenclature{$ \bm{n} = \left\{ n_{1}, n_{2}, n_{3} \right\}^{\operatorname{T}}$}{pseudo normal}
%
%
%
\begin{equation}
\bm{r}\left( \xi ,\eta ,\chi, t \right) = \bm{r}_m \left( \xi ,\eta, t \right) 
+ \chi \bm{n} \left( \xi ,\eta, t \right),  \quad \chi  \in \left[- d/2, d/2 \right], 
\label{eq: position}
\end{equation}
where $d$ is the thickness, 
$ \bm{r}_m $ the position of the material point at the mid-surface 
with the subscript $\left( \bullet \right)_m$ denoting the mid-surface.
Note that since the thickness is assumed to be constant during deformation
and the pseudo normal $\bm n$ represents the plate thickness direction, 
the distance $\chi$ is always between $- d/2$ and $d/2$. 
Therefore, the displacement $ \bm{u}$ of the material point can be determined by 
\begin{equation}
\bm{u} \left( \xi ,\eta ,\chi, t \right) = \bm{u}_m \left( \xi ,\eta, t \right) + \chi \Delta \bm{n} \left( \xi ,\eta, t  \right), 
\end{equation}
where $\Delta \bm{n} = \bm{n} - \bm{n}^0 $. 
\nomenclature{$\left( \bullet \right)_m$}
{representing the mid-surface parameter}
%
%
%
Then we can define the deformation gradient tensor as 
\begin{equation}\label{deformation_tensor}
\mathbb{F}  =  \nabla^0 \bm{r} = \nabla^0 \bm{u}  + \mathbb{I} 
= \left( \bm{a}_1, \bm{a}_2, \bm{a}_3 \right),
\end{equation}
where $\nabla^0 \equiv \partial / \partial \bm{\xi}$ 
is the gradient operator with respect to the initial configuration, 
$\mathbb{I}$ the identity matrix, 
and $\bm{a}_1 $,  $\bm{a}_2$,  $\bm{a}_3$ are specified by
\begin{equation}
\begin{cases}
\bm{a}_1  = \bm{r}_{m,\xi}  + \chi \bm{n}_{\xi} \\ 
\bm{a}_2  = \bm{r}_{m,\eta}  + \chi \bm{n}_{\eta}\\ 
\bm{a}_3  = \bm{n}
\end{cases}
\end{equation}
with 
$\nabla^0 \bm{r}_{m} \equiv (\bm{r}_{m,\xi}, \bm{r}_{m,\eta})^{\operatorname{T}}$
and 
$\nabla^0 \bm{n} \equiv (\bm{n}_{\xi}, \bm{r}_{\eta})^{\operatorname{T}}$. 
The deformation gradient tensor can be decomposed into two components as
\begin{equation}
	\mathbb{F}  =  \mathbb{F}_m + \chi \mathbb{F}_n,
\end{equation}
where $\mathbb{F}_m =  \left( \bm{r}_{m,\xi}^{\operatorname{T}}, \bm{r}_{m,\eta}^{\operatorname{T}}, \bm{n}^{\operatorname{T}} \right)$
and  $\mathbb{F}_n =  \left( \bm{n}_{\xi}^{\operatorname{T}},  \bm{n}_{\eta}^{\operatorname{T}}, 0 \right)$.
Furthermore, 
the real normal $\bm{n}_r$ is given as 
\begin{equation}
	\bm{n}_r = \frac{\bm{r}_{m,\xi} \times \bm{r}_{m,\eta}}{\left|\bm{r}_{m,\xi} \times \bm{r}_{m,\eta}\right|}.
\end{equation}

\subsubsection{Constitutive relation}
\label{sec:constitutive_relation}
With the deformation gradient tensor $\mathbb{F}$, 
the Green-Lagrangian strain tensor $\mathbb{E}$ can be obtained as 
\begin{equation}\label{Lagrangian-strain}
\mathbb{E}  = \frac{1}{2} \left(\mathbb{F}^{\operatorname{T}}\mathbb{F} - \mathbb{I}\right) 
= \frac{1}{2} \left(\mathbb{C} - \mathbb{I}\right),
\end{equation}
where $\mathbb{C}$ is the right Cauchy deformation gradient tensor. 
The Eulerian Almansi strain \fancy{$\epsilon$} 
can be converted from $\mathbb{E}$ as 
\begin{equation}\label{Almansi_strain}
	\fancy{$\epsilon$} = \mathbb{F}^{-\operatorname{T}} \cdot \mathbb{E} \cdot \mathbb{F}^{-1}
	= \frac{1}{2} \left(\mathbb{I} - \mathbb{F}^{-\operatorname{T}} \mathbb{F}^{-1}\right).
\end{equation}
When the material is linear and isotropic, 
the Cauchy stress \fancy{$\sigma$} reads 
\begin{equation}\label{constitutive_relation}
	\begin{split}
	   \fancy{$\sigma$} &= K \operatorname{tr}\left(\fancy{$\epsilon$} \right) \mathbb{I} + 2G\left( \fancy{$\epsilon$} - \frac{1}{3} \operatorname{tr}\left(\fancy{$\epsilon$} \right) \mathbb{I} \right) \\ 
		&= \lambda\operatorname{tr}\left( \fancy{$\epsilon$} \right) \mathbb{I} + 2\mu \fancy{$\epsilon$}, \\ 
	\end{split}
\end{equation}
where $\lambda$ and $\mu$ are the Lamé constants, 
$K = \lambda  + 2\mu /3$ the bulk modulus 
and $G = \mu$ the shear modulus. 
The relationship between the two moduli is given by 
\begin{equation}
	E = 2G \left( 1 + \nu \right) = 3K\left( 1 - 2\nu \right),
\end{equation}
where $E$ denotes the Young's modulus and $\nu$ the Poisson's ratio.
\subsubsection{Stress correction}\label{sec:stress_correction}
As the thickness is significantly 
less than the length and width of plate, 
the following boundary conditions hold when the plate 
is free from external forces on its surfaces 
where $\chi  =  \pm \frac{d}{2}$ or $z  =  \pm \frac{d}{2}$
\begin{equation}\label{shear_stress_correction}
	\left. \sigma^l_{xz}\right|_{z  =  \pm \frac{d}{2}} = 0, 
	\quad \left. \sigma^l_{yz}\right|_{z  =  \pm \frac{d}{2}} = 0,
\end{equation}
\begin{equation}\label{plane_stress}
	\left. \sigma^l_{zz}\right|_{z  \in \left[- \frac{d}{2} , \frac{d}{2}\right]} = 0,
\end{equation}
with the superscript $\left( \bullet \right)^l$ 
denoting the current local coordinates. 
Taking the boundary condition Eq. \eqref{plane_stress} 
and constitutive Eq. \eqref{constitutive_relation} into account,
the following relation of strains holds \cite{donning1998meshless}
\begin{equation}\label{strain_relation}
	\bar \epsilon^l_{zz}=\frac{-\nu \left( \epsilon^l_{xx} + \epsilon^l_{yy} \right)}{1-\nu},
\end{equation}
where the current local strain $\fancy{$\epsilon$}^l$ is obtained by 
\begin{equation}
	\fancy{$\epsilon$}^l = \mathbb{Q} \fancy{$\epsilon$} \mathbb{Q}^{\operatorname{T}}.
\end{equation}
Here, $\mathbb{Q}$ is the orthogonal transformation matrix 
from the global to current local coordinates. 
Following Batoz and Dhatt \cite{batoz1990modelisation},
$\mathbb{Q}$ can be given as
\begin{equation}
	\label{3D_transformation_matrix}
	\mathbb{Q} = \begin{bmatrix}
		n_3  + \frac{(n_2)^2 }{1 + n_3 } & - \frac{n_1 n_2}{1 + n_3 } & - n_1 \\
		- \frac{n_1 n_2 }{1 + n_3 } & n_3 + \frac{(n_1)^2 }{1 + n_3 }& -n_2 \\
		n_1 & n_2 & n_3 \\
	\end{bmatrix}.
\end{equation}
To satisfy the boundary conditions of Eq. \eqref{shear_stress_correction},
the transverse shear stress should be corrected as \cite{wisniewski2010finite}
\begin{equation}\label{shear_correction}
	\bar\sigma^l_{xz} = \bar\sigma^l_{zx} = \kappa \sigma^l_{xz},  
	\quad  \bar\sigma^l_{yz} = \bar\sigma^l_{zy}  =  \kappa \sigma^l_{yz},
\end{equation}
where $\kappa$ denotes the shear correction factor 
which is typically set to $5/6$ for the rectangular section of the isotropic plate.
Taking the corrected strain $\bar{\fancy{$\epsilon$}}^l$ into constitutive Eq. \eqref{constitutive_relation}
and then applying Eq. \eqref{shear_correction},
the corrected current local Cauchy stress $\bar{\fancy{$\sigma$}}^l$ is obtained.
%
%
%
%
%
%
%
%
%
%
%
%
%
%
%

\subsubsection{Conservation equations}
\label{}
The mass conservation equation can be written as 
\begin{equation}\label{mass-conservation}
	\rho  =  J_m^{-1}\rho^0,
\end{equation}
where $J_m = \det(\mathbb{F}_m)$,
$\rho^0$ and $\rho$ represent the initial and current densities, respectively. 
%
%
The momentum conservation equation is 
\begin{equation}\label{momentum-conservation1}
	\rho \ddot {\bm{u}}^l = \nabla  \cdot \left(\bar{ \fancy{$\sigma$}}^l\right)^{\operatorname{T}}
\end{equation}
or
\begin{equation}\label{momentum-conservation2}
	\rho \begin{bmatrix}
		\ddot u^l \\ \ddot v^l \\ \ddot w^l
	\end{bmatrix}
	= \begin{bmatrix}
		\frac{\partial \bar \sigma_{xx}^l}{\partial x} + \frac{\partial \bar \sigma_{xy}^l}{\partial y} 
		+ \frac{\partial \bar \sigma_{xz}^l}{\partial z} \\ 
		\frac{\partial \bar \sigma_{y x}^l}{\partial x} + \frac{\partial \bar \sigma_{yy}^l}{\partial y} 
		+ \frac{\partial \bar \sigma_{yz}^l}{\partial z}  \\ 
		\frac{\partial \bar \sigma_{z x}^l}{\partial x} + \frac{\partial \bar \sigma_{zy}^l}{\partial y} + \frac{\partial \bar \sigma_{zz}^l}{\partial z} 
	\end{bmatrix}.
\end{equation}
%
%
%
With Eqs. \eqref{shear_stress_correction} and \eqref{plane_stress}, 
we can integrate Eq. \eqref{momentum-conservation2} 
along $\chi$ or $z  \in \left[- d/2, d/2 \right]$ as 
\begin{equation}\label{momentum-conservation-integrated}
	d \rho  \begin{bmatrix}
		\ddot u_m^l \\ \ddot v_m^l \\ \ddot w_m^l
	\end{bmatrix}
	= \begin{bmatrix}
		\frac{\partial N_{xx}^l}{\partial x} + \frac{\partial N_{xy}^l}{\partial y} \\ 
		\frac{\partial N_{yx}^l}{\partial x} + \frac{\partial N_{yy}^l}{\partial y}\\ 
		\frac{\partial N_{zx}^l}{\partial x} + \frac{\partial N_{zy}^l}{\partial y}
	\end{bmatrix}, 
\end{equation}
where the stress resultant $\mathbb{N}^l$ is 
calculated by the Gauss–Legendre quadrature rule as 
\begin{equation} \label{gaussian_quadrature_rule1}
		\mathbb{N}^l  = \int_{-d/2}^{d/2} \bar {\fancy{$\sigma$}}^l \left(z\right) dz 
		= \sum\limits_{ip = 1}^N \bar {\fancy{$\sigma$}}^l \left(z_{ip}\right) A_{ip}.
\end{equation}
Here,  $z_{ip}$ is the integral point, 
$A_{ip}$ the weight, 
and $N$ the number of the integral point. 
Since the quadrature rule is conducted to yield an exact result 
for polynomials of degree $2 N - 1$ or lower \cite{gil2007numerical}, 
$N$ is determined by the applied constitutive relation.

By multiplying both sides of Eq. \eqref{momentum-conservation1} by $z$ 
and integrating along $z  \in \left[- d/2, d/2 \right]$, 
the angular momentum conservation equation can be obtained as
\begin{equation}\label{angular_momentum_conservation_integrated}
	\frac{d^3}{12} \rho \begin{bmatrix}
		\ddot n_1^l \\ 
		\ddot n_2^l \\
		\ddot n_3^l
	\end{bmatrix} 
	= \begin{bmatrix}
		\frac{\partial M_{xx}^l}{\partial x} + \frac{\partial M_{xy}^l}{\partial y} \\
		\frac{\partial M_{yx}^l}{\partial x} + \frac{\partial M_{yy}^l}{\partial y} \\ 
		\frac{\partial M_{zx}^l}{\partial x} + \frac{\partial M_{zy}^l}{\partial y} \\ 
	\end{bmatrix} 
	+ \begin{bmatrix}
		-N_{xz}^l \\ -N_{yz}^l \\ 0 
	\end{bmatrix},
\end{equation}
where the moment resultant $\mathbb{M}^l$ is calculated as 
\begin{equation} \label{gaussian_quadrature_rule2}
	\mathbb{M}^l  = \int_{-d/2}^{d/2} z \bar {\fancy{$\sigma$}}^l \left(z\right) dz 
	= \sum\limits_{ip = 1}^N z_{ip} \bar {\fancy{$\sigma$}}^l \left(z_{ip}\right) A_{ip}.
\end{equation}
%
%
%
Note that 
\begin{equation} 
	\int_{-d/2}^{d/2} z \frac{\partial \bar \sigma_{xz}^l}{\partial z} dz 
	= {\Big [}z \bar \sigma_{xz}^l {\Big ]}_{-d/2}^{d/2} - \int_{-d/2}^{d/2} \bar \sigma_{xz}^l dz = -N_{xz}^l.
\end{equation}

Therefore, the two governing equations, including 
the evolution of mid-surface displacement and pseudo normal, 
respectively,
for the 3D plate can be described as
\begin{equation}
\begin{cases}
	d\rho \bm{\ddot u}_m^l  = \nabla^l  \cdot \left(\mathbb{N}^l\right) ^{\operatorname{T}} \\
	\frac{d^3}{12}\rho \bm{\ddot n}^l  = \nabla^l \cdot \left(\mathbb{M}^l\right)^{\operatorname{T}}+ {\bf{Q}}^l,
\end{cases}
\end{equation}
where 
\begin{equation} \label{resultants}
	\mathbb{N}^l  = \begin{bmatrix}
		N_{xx}^l & N_{xy}^l & 0   \\
		N_{yx}^l & N_{yy}^l & 0   \\
		N_{zx}^l & N_{zy}^l & 0  \\
	\end{bmatrix},
	\mathbb{M}^l  = \begin{bmatrix}
		M_{xx}^l & M_{xy}^l & 0  \\
		M_{yx}^l & M_{yy}^l & 0  \\
		M_{zx}^l & M_{zy}^l & 0  \\
	\end{bmatrix},
	{\bf{Q}}^l  = \begin{bmatrix}
		-N_{xz}^l  \\ 
		-N_{yz}^l  \\ 
		0 \\ 
	\end{bmatrix}.
\end{equation}
In total Lagrangian formulation, 
the conservation equations above are converted into
\begin{equation} \label{plate-conservation-equation}
	\begin{cases}
		d\rho ^0 \bm{\ddot u}_m  
		= \left(\mathbb{F}_m\right)^{\operatorname{-T}} \nabla^0 
		\cdot \left( J_m \mathbb{N}^{\operatorname{T}} \right) \\
		\frac{d^3}{12}\rho^0 \bm{\ddot n}  
		= \left(\mathbb{F}_m\right)^{\operatorname{-T}} \nabla^0  
		\cdot \left( J_m \mathbb{M}^{\operatorname{T}} 		 \right)
		+ J_m \mathbb{Q}^{\operatorname{T}}{\bf{Q}}^l, 
	\end{cases}
\end{equation}
where $\mathbb{N} = \mathbb{Q}^{\operatorname{T}} \mathbb{N}^l \mathbb{Q}$ 
and $\mathbb{M} = \mathbb{Q}^{\operatorname{T}} \mathbb{M}^l \mathbb{Q}$ 
are the stress and moment resultants, respectively, in global coordinates. 
\subsection{3D shell model}
Based on the 3D plate model, 
the 3D shell model is obtained by introducing
the initial local coordinate system and the transformation matrix 
from the global to initial local coordinate system.
\subsubsection{Kinematics}
\begin{figure}[htb!]
	\centering
	\includegraphics[trim = 2mm 8mm 2mm 2mm, width = 0.9\textwidth]{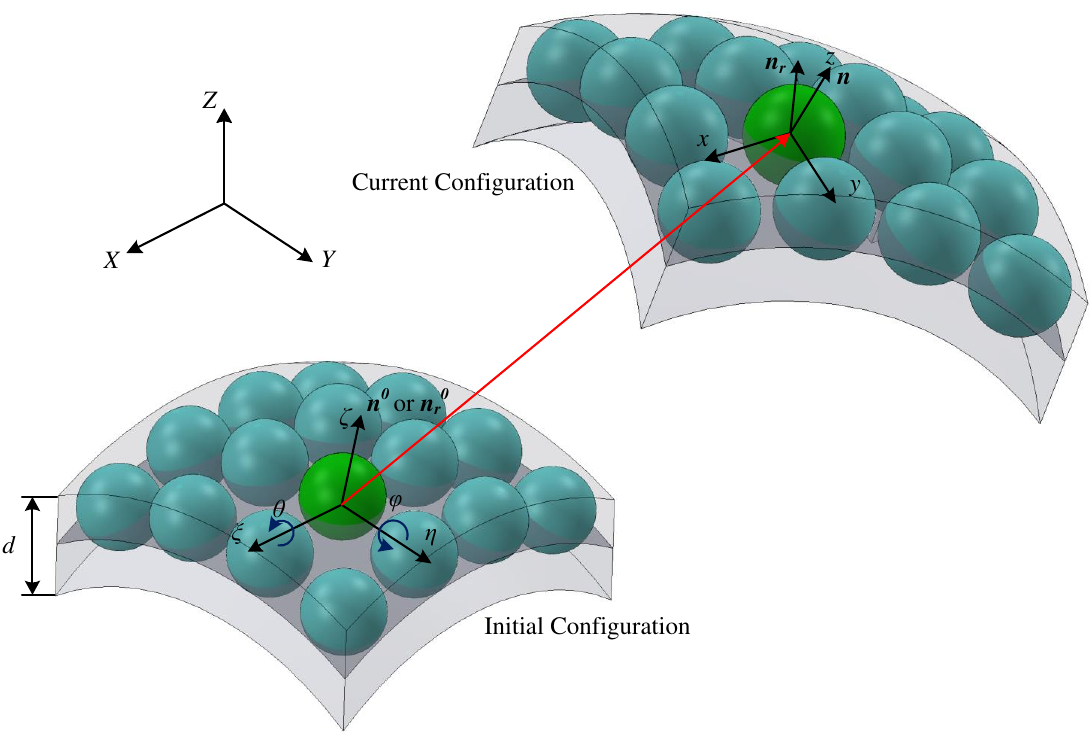}
	\caption{Schematic of a 3D shell model.}
	\label{figs:discretization_shell}
\end{figure}

The kinematics of shell can be constructed in the initial local coordinates 
denoted with the superscript $\left( \bullet \right)^L$.
Each material point possesses five degrees of freedom, viz.,
three translations $ \bm{u}^L  = \left\{ u^L, v^L, w^L \right\}^{\operatorname{T}} $ 
and two rotations
$\bm{\theta}^L  = \left\{ \theta^L, \varphi^L\right\}^{\operatorname{T}} $
as shown in Figure \ref{figs:discretization_shell}.
The pseudo-normal vector is also presented in initial local coordinates by 
$ \bm{n}^L = \left\{ n^L_1,n^L_2,n^L_3 \right\}^{\operatorname{T}}$,
especially denoted by $ \bm{n}^{0, L} = \left\{ 0, 0, 1 \right\}^{\operatorname{T}} $ 
in the initial local configuration.
The local position $ \bm{r}^L $ of a material point can be expressed as 
\begin{equation}
	\bm{r}^L\left( \xi ,\eta ,\chi, t \right) 
	= \bm{r}^L_m \left( \xi ,\eta, t \right) 
	+ \chi \bm{n}^L \left( \xi ,\eta, t \right),  
	\quad \chi  \in \left[- d/2,d/2 \right].
\end{equation}
The local displacement $ \bm{u}^L$ can thus be obtained by 
\begin{equation}
	\bm{u}^L \left( \xi ,\eta ,\chi, t \right) 
	= \bm{u}^L_m \left( \xi ,\eta, t \right) + \chi \Delta \bm{n}^L \left( \xi ,\eta, t  \right),
\end{equation}
where $\Delta \bm{n}^L = \bm{n}^L - \bm{n}^{0, L} $.
Similar to 3D plates, 
the local deformation gradient tensor of 3D shells can be defined as
\begin{equation}\label{shell_deformation_tensor}
	\mathbb{F}^L  = \nabla^{0, L} \bm{r}^L  + \nabla^{0, L} \bm{n}^L
	- \nabla^{0, L} \bm{n}^{0, L}
	= \left( {\bm{a}^L_1}, {\bm{a}^L_2}, {\bm{a}^L_3} \right),
\end{equation}
where $\nabla^{0, L} \equiv  \partial / \partial \bm{\xi}$ 
is the gradient operators defined in the initial local configuration, 
and $\bm{a}^L_1 $,  $\bm{a}^L_2$,  $\bm{a}^L_3$ are detailed by 
\begin{equation}
	\begin{cases}
		\bm{a}^L_1  = \bm{r}^L_{m,\xi}  + \chi \bm{n}^L_{\xi} - \chi \bm{n}^{0, L}_{\xi}  \\ 
		\bm{a}^L_2  = \bm{r}^L_{m,\eta}  + \chi \bm{n}^L_{\eta} - \chi \bm{n}^{0, L}_{\eta}  \\ 
		\bm{a}^L_3  = \bm{n}^L.
	\end{cases}
\end{equation}

\subsubsection{Stress correction and conservation equation}
With the local deformation gradient tensor $\mathbb{F}^L$, 
the local Eulerian Almansi strain $\fancy{$\epsilon$}^L$ 
can be calculated by the Eq. \eqref{Almansi_strain}.
After that,
the current local $\fancy{$\epsilon$}^l$ is obtained 
according to the coordinate transformation as
\begin{equation}
	\fancy{$\epsilon$}^l = \mathbb{Q} \left(\mathbb{Q}^0\right)^{\operatorname{T}} \fancy{$\epsilon$}^L  
	\mathbb{Q}^0 \mathbb{Q}^{\operatorname{T}},
\end{equation}
where $\mathbb{Q}^0$, 
the orthogonal transformation matrix from the global to initial local coordinates, 
is calculated from Eq. \eqref{3D_transformation_matrix} 
while the current pseudo normal $\bm{n}$ is replaced by the initial one $\bm{n}^0$. 
And then the corrected strain
$\bar{\fancy{$\epsilon$}}^l$ 
is estimated by applying Eq. \eqref{strain_relation}. 
After getting the current local Cauchy stress $\fancy{$\sigma$}^l$ 
by Eq. \eqref{constitutive_relation}, 
the corrected one $\bar{\fancy{$\sigma$}}^l$ is obtained by Eq. \eqref{shear_correction}.

Note that the total Lagrangian conservation equations of a 3D shell
has the same form as Eqs. \eqref{plate-conservation-equation}
with $\mathbb{F}_m = \left(\mathbb{Q}^0\right)^{\operatorname{T}} \mathbb{F}_m^L \mathbb{Q}^0$.
\subsection{2D plate/shell model}
If a plate/shell is assumed to be a strip 
that is very long and has a finite width,
and the transverse load is assumed to be uniform along the length, 
the analysis can be simplified 
at any cross section as a plane strain problem \cite{reddy2006theory}. 

The kinematics of 2D plate and shell 
can all be built in initial local coordinates, 
as the transformation matrix $\mathbb{Q}^0$ 
from the global to initial local coordinates
for the plate is the unit matrix.
The 2D model is in the global $X$-$Z$ plane, 
and each material point possesses three degrees of freedom, viz.,
two translations $ \bm{u}^L  = \left\{ u^L, w^L \right\}^{\operatorname{T}} $ 
and one rotation
$\bm{\theta}^L  = \left\{\varphi^L \right\}^{\operatorname{T}} $
expressed in the initial local coordinates.
The pseudo-normal vector is presented in the initial local coordinates 
by $ \bm{n}^L = \left\{ n^L_1, n^L_3 \right\}^{\operatorname{T}}$,
especially denoted by $ \bm{n}^{0, L} = \left\{ 0, 1 \right\}^{\operatorname{T}} $ 
in the initial local configuration. 
Similar to 3D model, 
the local position $ \bm{r}^L$ of a material point can be expressed as 
\begin{equation}
	\bm{r}^L\left( \xi ,\chi, t \right) = \bm{r}^L_m \left( \xi, t \right) + \chi \bm{n}^L \left( \xi, t \right),  
	\quad \chi  \in \left[- d/2, d/2 \right], 
\label{eq: shell_position}
\end{equation}
the local displacement $ \bm{u}^L$ can be evaluated as 
\begin{equation}
	\bm{u}^L \left( \xi ,\chi, t \right) = \bm{u}^L_m \left( \xi, t \right) + \chi \Delta \bm{n}^L \left( \xi, t  \right).
\end{equation}
and the local deformation gradient tensor is written as 
\begin{equation}\label{2D_shell_deformation_tensor}
	\mathbb{F}^L 
	= \nabla^{0, L} \bm{r}^L  + \nabla^{0, L} \bm{n}^L
	- \nabla^{0, L} \bm{n}^{0, L}
	= \left( {\bm{a}^L_1}, {\bm{a}^L_3} \right),
\end{equation}
where $\bm{a}^L_1 $ and $\bm{a}^L_3$ are given by
\begin{equation}
	\begin{cases}
		\bm{a}^L_1  = \bm{r}^L_{m,\xi}  + \zeta \bm{n}^L_{\xi} - \zeta \bm{n}^{0, L}_{\xi}  \\
		\bm{a}^L_3  = \bm{n}^L.
	\end{cases}
\end{equation}
The coordinate transformation matrix $\mathbb{Q}$ 
from global to current local coordinates 
is simplified from Eqs. \eqref{3D_transformation_matrix} as
\begin{equation} \label{2D_Q}
	\mathbb{Q} = \begin{bmatrix}
		n_3 & -n_1 \\
		n_1 & n_3 \\
	\end{bmatrix},
\end{equation}
and the 2D transformation matrix $\mathbb{Q}^0$ from global to initial local coordinates 
can also calculated by Eq. \eqref{2D_Q} 
while the current pseudo normal $\bm{n}$ 
is replaced by the initial one $\bm{n}^0$.
The corrected relation of strains is simplified from Eq. \eqref{strain_relation} as
\begin{equation}\label{2D_strain_relation}
	\bar\epsilon^l_{zz}=\frac{-\nu  \epsilon^l_{xx} }
	{1-\nu}.
\end{equation}
Finally, the 2D conservation equation is identical to 3D  
Eq. \eqref{plate-conservation-equation} with
\begin{equation}
	\mathbb{N}^l  = \begin{bmatrix}
		N^l_{xx} & 0   \\
		N^l_{z x}& 0  \\
	\end{bmatrix},
	\mathbb{M}^l  = \begin{bmatrix}
		M^l_{xx} & 0  \\
		M^l_{zx}  & 0  \\
	\end{bmatrix},
	{\bf{Q}}^l  = \begin{bmatrix}
		-N^l_{xz}  \\ 
		0 \\ 
	\end{bmatrix}.
\end{equation}
\section{SPH method for plate and shell structures}\label{sec:SPH_method}
In this section, we first introduce 
the reduced-dimensional SPH method, 
and detail the proposed formulations for plate and shell structures,  
including the discretization of conservation equations, 
non-singular conversion algorithm for 
the kinematics between rotation angles and pseudo normal, 
and the algorithms to increase numerical stability and alleviate hourglass modes.
After that, the time-integration schemes are presented.
\subsection{Reduced-dimensional SPH method}
In full-dimensional SPH method, 
the smoothed field $f(\bm{r})$ is obtained as
\begin{equation} \label{eq:kenelintegral}
	f(\bm{r})  = \int_{\Omega} f(\bm{r}')  W(\bm{r} - \bm{r}', h) d\bm{r}',
\end{equation}
where $f(\bm{r}')$ is the original continuous field before smoothing,  
$\Omega$ the entire space and
$W(\bm{r} - \bm{r}', h)$ a Gaussian-like kernel function 
with smoothing length $h$ denoting the compact support. 
By carrying out the integration of Eq. \eqref{eq:kenelintegral} 
along the thickness of the plate/shell structure, 
we can obtained the reduced-dimensional smoothed field by
\begin{equation}\label{eq:reducedkenelintegral}
f \left( \mathbf r\right)   \approx \int_{\widehat \Omega}  f \left( \mathbf r' \right) \widehat  W \left( \mathbf r-\mathbf r', h\right) d\mathbf r', 
\end{equation}
where $\widehat \Omega$ denotes the reduced space 
and $\widehat W(\bm{r} - \bm{r}', h)$ the reduced kernel function. 
Note that Eqs. (\ref{eq:kenelintegral}) and (\ref{eq:reducedkenelintegral}) 
have identical forms of formulation. 
A reduced-dimensional fifth-order Wendland kernel 
\cite{wendland1995piecewise} reads
\begin{equation}
\widehat W (q, h) =  \alpha
\begin{cases}
\left(1 + 2q \right) \left( 1 - q / 2\right)^4  & \text{if} \quad 0\le q\le2 \\ 
0 & \text{otherwise} 
\end{cases},
\end{equation}
where $q = \left| \mathbf r-\mathbf r'\right|/h$ and the constant $\alpha$ is equal to $\frac{3}{4 h}$ and $\frac{7}{4\pi h^2}$ 
for 2D and 3D problems, respectively. 
Also note that the reduced kernel function has identical form with the full-dimensional counterpart except 
different dimensional normalizing constant parameter, 
allowing the integration of unit can be satisfied 
in the reduced space.
Due to the almost identical forms, 
in present work from here, 
we do not identify the full- and reduced-dimensional formulations 
unless explicitly mentioned. 

In the reduced-dimensional SPH method, 
similarly to the full-dimensional counterpart  \cite{monaghan2005smoothed}, 
the gradient of the variable field $f(\bm{r})$ 
at a surface particle $i$ can be approximated as 
\begin{equation}
	\label{eq:gradsph}
	\begin{split}
		\nabla f_i & = \int_{\Omega} \nabla f (\bm{r}) W(\bm{r}_i - \bm{r}, h) d \bm{r}  \\
		& =  - \int_{\Omega} f (\bm{r}) \nabla W(\bm{r}_i - \bm{r}, h) d \bm{r}
		\approx  - \sum\limits_j  f_j\nabla W_{ij}V_j ,
	\end{split}
\end{equation}
where $V$ is the reduced particle volume, i.e. length 
and area for 2D and 3D problems, respectively. 
Here, the summation is conducted over all the neighboring particles $j$ 
located at the support domain of the particle $i$, 
and $\nabla W_{ij}  = -\frac{\partial W\left( \bm{r}_{ij}, h \right)}{\partial r_{ij}} \bm{e}_{ij} $  
is the gradient of the kernel function with 
$\bm{r}_{ij} = \bm{r}_{i} - \bm{r}_{j}$
and $\bm{e}_{ij} = \bm{r}_{ij} / |\bm{r}_{ij}|$ denoting 
the unit vector pointing from particle $j$ to $i$.
Equation \eqref{eq:gradsph} can be modified into a strong form as
\begin{equation}
	\label{eq:gradsph-strong}
	\nabla f_i = \nabla f_i - f_{i}\nabla 1  \approx   \sum\limits_j  f_{ij} \nabla W_{ij}V_j, 
\end{equation}
where $f_{ij} = f_{i} - f_{j}$ is the interparticle difference value. 
This strong-form derivative operator can be used to determine the local structure of a field, such as the deformation gradient tensor. 
And Eq. \eqref{eq:gradsph} can also be modified into a weak form as 
\begin{equation}
	\label{eq:gradsph-weak}
	\nabla f_i = f_{i}\nabla 1 + \nabla f_i \approx  - \sum\limits_j \left(f_{i} + f_{j}\right)\nabla W_{ij} V_j. 
\end{equation}
This weak-form derivative operator is applied 
here for solving the conservation equations. 
Thanks to its anti-symmetric feature, 
i.e., $\nabla W_{ij} = - \nabla W_{ji}$, 
the momentum conservation of the particle system 
is ensured \cite{monaghan2005smoothed}. 

Also note that, 
in the present work, 
the reduced-dimensional SPH method is 
employed for total Lagrangian formulation \cite{randles1996smoothed},
such as Eq. \eqref{plate-conservation-equation}. 
Therefore, 
the smoothing kernel function and 
its derivatives are only evaluated once, 
also denoted with superscript $\left( \bullet \right)^0$ 
at the initial configuration, 
and kept unchanged during the simulation.  
\subsection{First-order consistency corrections}
\label{sec:consistency correction}
For the full-dimensional SPH in total Lagrangian formulation,
in order to remedy the 1st-order inconsistency 
which is caused by incomplete kernel support 
at domain boundary 
or with irregular particle distribution, 
the symmetric correction matrix $\mathbb{B}^{0}_i$ for each particle
\cite{vignjevic2006sph, liu2010smoothed}
is introduced for each particle to satisfy the linear-reproducing condition
\begin{equation}\label{eq:consistency-condition}
\left(\sum\limits_j {\bm{r}^{0}_{ij} \otimes \nabla^{0} W_{ij} V^{0}_j} \right)
\mathbb{B}^{0}_i  = \mathbb{I}.
\end{equation}
Then the strong-form approximations of gradient Eq. \eqref{eq:gradsph-strong} is
modified as
\begin{equation}
\label{eq:gradsph-strong-corrected}
\nabla^{0} f_i \approx \left( \sum\limits_j  f_{ij} \nabla^{0} W_{ij}V^{0}_j \right) \mathbb{B}^{0}_i, 
\end{equation}
and the weak-form approximations of divergence Eq.  \eqref{eq:gradsph-weak} as
\begin{equation}
\label{eq:gradsph-weak-corrected}
\nabla^{0} \cdot f_i \approx  - \sum\limits_j \left(f_{i}\mathbb{B}^{0}_i + f_{j}\mathbb{B}^{0}_j\right)\nabla^{0} W_{ij} V^{0}_j.
\end{equation}
In  the reduced-dimensional SPH,
we generalize the linear-reproducing condition as
\begin{equation}\label{eq:reduced-consistency-condition}
\left[\mathbb{G}^{\operatorname{T}} \mathbb{Q}^{0}_i \left(\sum\limits_j {\bm{q}^{0}_{ij} \otimes \nabla^{0} W_{ij} V^{0}_j} \right) \left(\mathbb{Q}^{0}_i\right)^{\operatorname{T}} \mathbb{G} \right] \mathbb{B}^{0, L}_i  = \mathbb{K}_i,
\end{equation}
where $\bm{q}^{0}_{ij}$ is 
the initial inter-particle difference of a linear vector,
$\mathbb{Q}^{0}_i$ 
is the transformation matrix from the global to initial local coordinates, 
and $\mathbb{G}$ is a reducing matrix, i.e.,
\begin{equation}
\mathbb{G} = 
\begin{bmatrix}
1  \\
0  \\
\end{bmatrix}
~\text{and}~
\begin{bmatrix}
1 & 0  \\
0 & 1  \\
0 & 0  \\
\end{bmatrix} 
\end{equation}
for 2D and 3D problems, respectively.
It ensures that the corrections are carried out 
within the local reduced space.
Similarly, the strong-form approximations of gradient Eq. \eqref{eq:gradsph-strong} is
modified as
\begin{equation}
\label{eq:shell-strong-corrected}
\nabla^{0} f_i \approx  \left( \sum\limits_j  f_{ij} \nabla^{0} W_{ij}V^{0}_j \right)\widetilde{\mathbb{B}}^{0}_i, 
\end{equation}
where $\widetilde{\mathbb{B}}^{0}_i = \left(\mathbb{Q}^{0}_i\right)^{\operatorname{T}} \mathbb{G} \mathbb{B}^{0, L}_i \mathbb{G}^{\operatorname{T}} \mathbb{Q}^{0}_i$ 
and the weak-form approximations of divergence Eq.  \eqref{eq:gradsph-weak} as
\begin{equation}
\label{eq:shell-weak-corrected}
\nabla^{0} \cdot f_i \approx  - \sum\limits_j \left(f_{i}\widetilde{\mathbb{B}}^{0}_i + f_{j}\widetilde{\mathbb{B}}^{0}_j\right)\nabla^{0} W_{ij} V^{0}_j.
\end{equation}
Here, we introduce the correction matrix 
$\widetilde{\mathbb{B}}^{0}_i = \widetilde{\mathbb{B}}^{0,\bm{r}}_i$,
 $\bm{q}^{0}_{ij} = \bm{r}_{ij}^0$ 
and $\mathbb{K}_i$ is the reduced identity matrix denoted as
\begin{equation}
\mathbb{K}_i = \mathbb{K}^{\bm{r}} = 
\begin{bmatrix}
1 
\end{bmatrix}
~\text{and}~
\begin{bmatrix}
1 & 0  \\
0 & 1  
\end{bmatrix} 
\end{equation}
for 2D and 3D problems, respectively, 
to correct the position-based quantities.
Similarly, we introduce the correction matrix 
$\widetilde{\mathbb{B}}^{0}_i = \widetilde{\mathbb{B}}^{0,\bm{n}}_i$,
$\bm{q}^{0}_{ij} = \bm{n}_{ij}^0$ 
and 
\begin{equation}
\mathbb{K}_i = \mathbb{K}^{\bm{n}}_{i} = 
\begin{bmatrix}
1/R_i^L
\end{bmatrix}
~\text{and}~
\begin{bmatrix}
1/R_{1,i}^L & 0 \\
0 & 1/R_{2,i}^L 
\end{bmatrix},
\end{equation}
where $R_i^L$, $R_{1,i}^L$ and $R_{2,i}^L$ 
are the curvature radii of particle $i$
for 2D and 3D problems, respectively, 
to correct rotation-based quantities. 
\subsection{Discreization of conservation equations}
\label{sec:reduced_dimensional_TLSPH}
With two correction matrices obtained from 
Eq. \eqref{eq:reduced-consistency-condition} 
and following Eq. \eqref{eq:shell-weak-corrected}, 
the momentum equations \eqref{plate-conservation-equation} 
are discretized  as
\begin{equation}\label{discrete_dynamic_equation1}
	d\rho _{i}^0 \bm{\ddot u}_{m, i} =	\sum\limits_j 
	\left(
	J_{m, i} \mathbb{N}_i
	 \left(\mathbb{F}_{m, i}  \right)^{\operatorname{-T}} 
	 \widetilde{\mathbb{B}}^{0, \bm{r}}_i 
	+ J_{m, j} \mathbb{N}_j \left(\mathbb{F}_{m, j}  \right)^{\operatorname{-T}}  \widetilde{\mathbb{B}}^{0, \bm{r}}_j
	\right)
	\nabla^0 W_{ij} V_j^0,
\end{equation} 
and
\begin{equation}\label{discrete_dynamic_equation2}
	\begin{split}
\frac{d^3 } {12}	\rho_i^0\bm{\ddot n }_i &  = 
	\sum\limits_j 
\left(
	J_{m, i} \mathbb{M}_i
	\left(\mathbb{F}_{m, i}  \right)^{\operatorname{-T}} 
	\widetilde{\mathbb{B}}^{0, \bm{n}}_i + J_{m,j} \mathbb{M}_j
	\left(\mathbb{F}_{m, j}  \right)^{\operatorname{-T}} 
	\widetilde{\mathbb{B}}^{0, \bm{n}}_j
	\right)
	\nabla^0 W_{ij} V_j^0 \\
	& + {J_{m, i} \left(\mathbb{Q}_i^0\right)}^{\operatorname{T}} \bm{Q}^l_i. 
	\end{split}
\end{equation}

\subsection{Kelvin–Voigt type damping}
Following Ref. \cite{zhang2022artificial}, 
when calculating the current local Cauchy stress 
by using the constitutive Eq. \eqref{constitutive_relation},
an artificial damping stress $\fancy{$\sigma$}_d^l$ 
based on the Kelvin-Voigt type damper is introduced here as
%
%
%
\begin{equation}
	\label{Cauchy_stress_damping}
	\fancy{$\sigma$}_d^l = J_m^{-1} \mathbb{Q} \left(\mathbb{Q}^0\right)^{\operatorname{T}} 
		\mathbb{F}^L 
		\dot{\mathbb{E}}^L 
		\fancy{$\gamma$}
		\left(\mathbb{F}^L\right)^{\operatorname{T}} 
		\mathbb{Q}^0 \mathbb{Q}^{\operatorname{T}},
\end{equation}
where the numerical viscosity matrix 
\begin{equation}
	\fancy{$\gamma$} = 
	\begin{bmatrix}
		\rho c h / 2 & 0  \\
		0 & \rho c s / 2  \\
	\end{bmatrix}
	~\text{and}~
	\begin{bmatrix}
		\rho c h / 2 & 0 & 0  \\
		0 & \rho c h / 2 & 0  \\
		0 & 0 & \rho c s / 2  \\
	\end{bmatrix} 
\end{equation}
where $c  = \sqrt {K/\rho} $ and $s = \min (h, d)$,  
for 2D and 3D problems, respectively.
Note that, 
different from Ref. \cite{zhang2022artificial}, 
where an isotropic numerical damping is applied, 
the present damping leads to a smaller out-of-plane contribution 
when $d \le h$, 
which makes it suitable for 
both thin and moderately thick plate and shell structures.
The change rate of the Green-Lagrangian strain tensor is given as
\begin{equation}
	\dot{\mathbb{E}}^L
	= \frac{1}{2} \left[\left(\dot{\mathbb{F}}^L \right)^{\operatorname{T}} \mathbb{F}^L  
	+  \left(\mathbb{F}^L\right)^{\operatorname{T}}  \dot{\mathbb{F}}^L \right].
\end{equation}
Here, the change rate of the deformation gradient tensor of particle $i$ is
\begin{equation}\label{change_rate_deformation_gradient}
{\mathbb{\dot F}}^L_i = \nabla^{0, L} \bm{\dot u}_i^L 
= \nabla^{0}\bm{\dot u}^L_{m, i} + \chi \nabla^{0}\bm{\dot n}^L_{i} ,
\end{equation}
where
\begin{equation}
	\begin{cases}
		\nabla^{0}\bm{\dot u}^L_{m, i} = \mathbb{Q}_i^0\left(\sum\limits_j 
		\bm{\dot u}_{m,ij}
		\otimes \nabla^0 W_{ij} V_j^0  \right) \widetilde{\mathbb{B}}_i^{0, \bm{r}} \left(\mathbb{Q}_i^0\right)^{\operatorname{T}} \\
		\nabla^{0}\bm{\dot n}^L_{i} = \mathbb{Q}_i^0 \left(\sum\limits_j
		\bm{\dot n}_{ij} 
		\otimes \nabla^0 W_{ij} V_j^0 \right) \widetilde{\mathbb{B}}_{i}^{0,  \bm{n}} \left(\mathbb{Q}_i^0\right)^{\operatorname{T}}
	\end{cases}
\end{equation}
are obtained following the consistency condition Eq. \eqref{eq:reduced-consistency-condition} 
and the strong-form correction Eq. \eqref{eq:shell-strong-corrected}.
\subsection{Hourglass control}
Inspired from Refs. \cite{kondo2010suppressing, ganzenmuller2015hourglass} in full-dimensional SPH for total Lagrangian solid dynamics, 
we introduce a hourglass control algorithm here to 
alleviate the hourglass modes in the dynamics of plate and shell structures.

First, 
we estimate the position of the inter-particle middle point linearly
using the deformation gradient tensor from particles $i$ and $j$, 
respectively, as 
\begin{equation}\label{predicted_position1}
    \bm{r}_{i + \frac{1}{2}} =  \bm{r}_{m, i}
    - \frac{1}{2} \mathbb{F}_{m, i} \bm{r}_{m, ij}^0,
\quad
	\bm{r}_{j - \frac{1}{2}} =  \bm{r}_{m,j} 
	+ \frac{1}{2} \mathbb{F}_{m,j} \bm{r}_{m, ij}^0.
\end{equation}
One can find that the inconsistency beyond linear estimation 
$\hat{\bm{r}}_{ij} = \bm{r}_{i + \frac{1}{2}} - \bm{r}_{j - \frac{1}{2}}$ is
\begin{equation}
	\hat{\bm{r}}_{ij} = \bm{r}_{m, ij}
	- \frac{1}{2} \left( \mathbb{F}_{m,i} + \mathbb{F}_{m,j} \right) \bm{r}_{m, ij}^0.
\end{equation}
To suppress the position inconsistency  $\hat{\bm{r}}_{ij}$,
we introduce an extra correction term to the discrete momentum conservation Eq. \eqref{discrete_dynamic_equation1} as
\begin{equation}\label{artificial-stress}
		d \rho _{i}^0 \bm{\ddot u}^{cor}_{m,i} 
		= \sum\limits_j \alpha G \beta_{ij} \gamma_{ij}^{\bm{r}} D \hat{\bm{r}}_{ij} 
		\frac{\partial W\left( \bm{r}_{ij}^0, h \right)}
		{\partial r_{ij}^0} V_j^0 \\
\end{equation}
where $\beta_{ij} = W^0_{ij} / W_{0}$ 
leads to less correction to further neighbors, 
$\gamma_{ij}^{\bm{r}} 
= \min \left(2\left|\hat{\bm{r}}_{ij}\right|/
\left|\bm{r}_{m, ij}\right|, 1\right)$ is an adaptive limiter 
for less correction on the domain where the inconsistency is less pronounced, 
D the dimension, 
and parameter $\alpha = 0.002$ according to the numerical experiment 
and remains constant throughout this work. 
Note that, since the inconsistency decreases with decreasing particle spacing, 
different from Refs. \cite{ganzenmuller2015hourglass}, 
the present correction is purely numerical 
and vanishes with increasing resolution.
Similarly, for the predicted pseudo normal, 
the difference of the intermediate point can be described as 
\begin{equation}
	\hat{\bm{n}}_{ij} =  \bm{n}_{ij} - \bm{n}_{ij}^0
	- \frac{1}{2} \left(\mathbb{F}_{\bm{n},i} + \mathbb{F}_{\bm{n},j}\right) \bm{r}_{ij}^0.
\end{equation}
Similar with Eq. \eqref{artificial-stress},
the extra correction term added to 
the discrete angular momentum conservation Eq. \eqref{discrete_dynamic_equation2} is
\begin{equation}\label{artificial-torque}
	\frac{d^{3}}{12}\rho _{i}^0	
		\bm{\ddot n}^{cor}_i  
	 = \sum\limits_j \alpha  G d^2 \beta_{ij} \gamma_{ij}^{\bm{n}} D
	 \hat{\bm{n}}_{ij}
		\frac{\partial W\left( \bm{r}_{ij}^0,h \right)}{\partial r_{ij}^0 }
		V_j^0, 
\end{equation}
where the adaptive limiter is $\gamma_{ij}^{\bm{n}} 
= \min \left(2 \left|\hat{\bm{n}}_{ij}\right| /
\left|\bm{n}_{ij} - \bm{n}_{ij}^0\right|, 1\right)$.
Note that, 
different with Refs. \cite{kondo2010suppressing, ganzenmuller2015hourglass},
the present correction force is introduced in particle pairwise pattern, 
implying momentum conservation \cite{monaghan2005smoothed}.
Also note that, 
the correction force is scaled to the shear, 
rather than Young's, modulus, 
due to the fact that the hourglass modes are 
characterized by shear deformation \cite{wu2023essentially}.
\subsection{Conversion between rotations and pseudo normal} 
Different from the mid-surface displacement,
which can be numerically integrated directly from its evolution equation,
the pseudo normal is not suitable for direct numerical integration 
since its unit magnitude may not be maintained strictly.    
Under the assumption of small rotation, 
one may have the simplified relation between the pseudo normal and rotations,
i.e. $\bm{\ddot{\theta}} = (-\ddot{n}_2, \ddot{n}_1)$, 
so that one can obtain the rotation increment,  
and update rotation matrix $\mathbb{R}$ 
using Rodrigues formula 
\cite{betsch1998parametrization, lin2014efficient} and 
finally the integrated pseudo normal by $\bm{n} = \mathbb{R}\bm{n}^0$. 

In the present work, the numerical integration of pseudo normal is 
carried out without the assumption of small rotation,
that is, we strictly identify the rotations and pseudo normal 
by using the original evolution equation 
and obtain their conversion relations
$\bm{\ddot{\theta}} = \bm{\ddot{\theta}}(\bm{\ddot{n}}, \bm{\dot{\theta}}, \bm{\theta})$.
Different from using the rotation matrix $\mathbb{R}$ based on Rodrigues formula, 
we update the pseudo normal $\bm{n}^L$ 
with \cite{betsch1998parametrization, wisniewski2010finite}
\begin{equation} \label{eq:updated-pseudo-normal}
\bm{n}^L = \mathbb{R}^L_{\eta} \mathbb{R}^L_{\xi} \bm{n}^{0, L},
\end{equation}
where $\mathbb{R}_{\xi}^L\equiv \mathbb{R}_{\xi}(\theta^{L})$ and 
$\mathbb{R}^L_{\eta} \equiv \mathbb{R}_{\eta}(\varphi^{L})$ 
are the local rotation matrices
respected to the axes $\xi$ and $\eta$, respectively,
or, equivalently, with the change rate
\begin{equation} \label{eq:pseudo-normal-rate}
\bm{\dot{n}}^L = \dot{\mathbb{R}}^L_{\eta} \dot{\mathbb{R}}^L_{\xi},
\end{equation}
where 
$\dot{\mathbb{R}}^L_{\xi}\equiv \mathbb{R}_{\xi}(\theta^{L}, \dot{\theta}^L)$
and 
$\dot{\mathbb{R}}^L_{\eta}\equiv \mathbb{R}_{\eta}(\varphi^{L}, \dot{\varphi}^L)$.
Here, the rotations and their change rates are numerically integrated directly 
with the help of conversion relations.

Specifically, for a 2D problem,  
$\mathbb{R}^L_{\xi}$ is a unit matrix, 
and $\mathbb{R}^L_{\eta}$ can be described as
\begin{equation}
\mathbb{R}^L_{\eta}  =  \begin{bmatrix}
\cos \varphi^L & \sin \varphi^L   \\
- \sin \varphi^L & \cos \varphi^L   \\
\end{bmatrix}.
\end{equation}
Then, one has the relation as
\begin{equation} \label{2d-rotaion-pseudo-normal}
\bm{n}^L = (\sin \varphi^L, \cos \varphi^L)^{\operatorname{T}},
\end{equation}
its 1st-order time derivative corresponding Eq. \eqref{eq:pseudo-normal-rate}
\begin{equation}\label{2D_n_first_derivative}
\dot {\bm{n}}^L 
= (\cos \varphi^L \cdot \dot \varphi^L, 
-\sin \varphi^L \cdot \dot \varphi^L)^{\operatorname{T}}, 
\end{equation}
and 2nd-order derivative
\begin{equation}\label{2D_n_second_derivative}
\ddot {\bm{n}}^L 
= (-\sin \varphi^L \cdot \left(\dot \varphi^L\right)^2 
+\cos \varphi^L \cdot \ddot \varphi^L, 
-\cos \varphi^L \cdot \left(\dot \varphi^L\right)^2 
-\sin \varphi^L \cdot \ddot \varphi^L)^{\operatorname{T}}.
\end{equation}
Note that Eq. \eqref{2D_n_second_derivative} suggests 
two theoretically equivalent conversion relations 
\begin{equation}\label{2D_conversion_relation1}
\ddot \varphi^L 
= \frac{\ddot n_1^L + \sin \varphi^L \cdot \left(\dot \varphi^L\right)^2}
{\cos \varphi^L}
\quad \text{and} \quad
\ddot \varphi^L 
= \frac{\ddot n_2^L + \cos \varphi^L \cdot \left(\dot \varphi^L\right)^2}
{-\sin \varphi^L}.
\end{equation}
Although each of them can be used 
to obtain the rotation angle $\varphi^L$ and its change rate with 
direct numerical integration
and hence the pseudo normal with Eq. \eqref{eq:updated-pseudo-normal}, 
there are singularities at large rotation angles
$\varphi^L = 0.5 \pi + k \pi$ (1st relation) 
or $\varphi^L = k \pi$ (2nd relation) with $k=0, 1, 2, 3, ...$ 
\cite{simo1990stress, betsch1998parametrization, singla2004avoid}. 
In order to eliminate the singularities, 
we propose to uses both relations with a weighted average as  
\begin{equation}\label{2D_conversion_relation3}
\begin{split}
\ddot \varphi^L 
& = \left(\cos \varphi^L\right)^2 
\frac{\ddot n_1^L + \sin \varphi^L \cdot \left(\dot \varphi^L\right)^2}
{\cos \varphi^L}
+ \left(\sin \varphi^L\right)^2 
\frac{\ddot n_2^L + \cos \varphi^L \cdot \left(\dot \varphi^L\right)^2}
{-\sin \varphi^L} \\
& = \cos \varphi^L 
\left(\ddot n_1^L + \sin \varphi^L \cdot \left(\dot \varphi^L\right)^2\right)
- \sin \varphi^L 
\left(\ddot n_2^L + \cos \varphi^L \cdot \left(\dot \varphi^L\right)^2\right),
\end{split}
\end{equation}
which cancels both denominators. 
Note that the present formulation recovers 
$\ddot \varphi^L = \ddot n_1^L$ under the assumption of small rotation. 

As for 3D problems, 
the rotation matrices $\mathbb{R}^L_{\xi}$ and $\mathbb{R}^L_{\eta}$ are
\begin{equation}
\mathbb{R}^L_\xi  =  \begin{bmatrix}
1 & 0 & 0  \\
0 & \cos \theta^L  &  -\sin \theta^L  \\
0 & \sin \theta^L  & \cos \theta^L   \\
\end{bmatrix},
\end{equation}
and
\begin{equation}
\mathbb{R}^L_\eta  =  \begin{bmatrix}
\cos \varphi^L  & 0 & \sin \varphi^L   \\
0 & 1 & 0  \\
- \sin \varphi^L & 0 & \cos \varphi^L   \\
\end{bmatrix}.
\end{equation}
Similarly, 
one has the relation between rotations and pseudo normal \cite{hughes1978consistent}
\begin{equation}
\bm{n}^L = ( \cos \theta^L \sin \varphi^L, -\sin \theta^L,  \cos \theta^L \cos \varphi^L)^{\operatorname{T}},
\end{equation}
its 1st-order time derivatives corresponding Eq. \eqref{eq:pseudo-normal-rate}
\begin{equation}\label{3D_n_first_derivative}
\begin{cases}
\dot n^L_1 = -\sin \theta^L \sin \varphi^L \dot \theta^L 
+ \cos \theta^L \cos \varphi^L \dot \varphi^L\\
\dot n^L_2 = -\cos \theta^L \dot \theta^L\\
\dot n^L_3 = -\sin \theta^L \cos \varphi^L \dot \theta^L 
- \cos \theta^L \sin \varphi^L \dot \varphi^L,
\end{cases}
\end{equation}
and 2nd-order derivatives
\begin{equation}\label{3D_n_second_derivative}
\begin{cases}
\begin{split}
\ddot n^L_1 = 
& -\sin \theta^L \sin \varphi^L \ddot \theta^L
- \cos \theta^L \sin \varphi^L ({\dot\theta}^L)^2
- 2 \sin \theta^L \cos\varphi^L \dot \theta^L \dot \varphi^L \\
& - \cos \theta^L \sin \varphi^L (\dot \varphi^L)^2
+ \cos \theta^L \cos \varphi^L \ddot \varphi^L
\end{split}\\
\ddot n^L_2 = \sin \theta^L (\dot \theta^L)^2 -\cos \theta^L \ddot \theta^L\\
\begin{split}
\ddot n^L_3 
&=  -\sin \theta^L \cos \varphi^L \ddot \theta^L
- \cos \theta^L \cos \varphi^L ({\dot\theta}^L)^2
+ 2 \sin \theta^L \cos\varphi^L \dot \theta^L \dot \varphi^L \\
& - \cos \theta^L \cos \varphi^L (\dot \varphi^L)^2
- \cos \theta^L \sin \varphi^L \ddot \varphi^L.
\end{split}
\end{cases}
\end{equation}
Note that, one can obtained 3 theoretically equivalent conversion relations,
respectively,
by 1st and 3rd expressions of Eq. \eqref{3D_n_second_derivative} as
\begin{equation}\label{3D_n_second_derivative1}
\begin{cases}
\ddot \theta^L = 
-\left(
\ddot n^L_3 \cos \varphi^L+\ddot n^L_1 \sin \varphi^L 
+ \left(\dot \varphi^L \right)^2 \cos \theta^L + \left(\dot \theta^L \right)^2 \cos \theta^L
\right)  / \sin \theta^L\\
\ddot \varphi^L = 
\left(
\ddot n^L_1 \cos \varphi^L - \ddot n^L_3 \sin \varphi^L 
+ 2 \dot \varphi^L \dot \theta^L \sin \theta^L
\right)  / \cos \theta^L,
\end{cases}
\end{equation}
1st and 2nd expressions
\begin{equation}\label{3D_n_second_derivative2}
\begin{cases}
\ddot \theta^L = 
\left(
\sin \theta^L \left(\dot \theta^L \right)^2 - \ddot n^L_2
\right)  / \cos \theta^L\\
\begin{split}
\ddot \varphi^L 
&= (
\ddot n^L_1 \cos \theta^L 
+ \left(\dot \varphi^L\right)^2 \cos^2 \theta^L \sin \varphi^L
+ \left(\dot \theta^L\right)^2 \sin \varphi^L
- \ddot n^L_2 \sin \varphi^L \sin \theta^L \\
&+ 2 \dot \varphi^L \dot \theta^L \cos \varphi^L \cos\theta^L \sin \theta^L
)  / \cos \varphi^L \cos^2 \theta^L,
\end{split}
\end{cases}
\end{equation}
and 2nd and 3rd expressions
\begin{equation}\label{3D_n_second_derivative3}
\begin{cases}
\ddot \theta^L = 
\left(
\sin \theta^L \left(\dot \theta^L \right)^2 - \ddot n^L_2
\right)  / \cos \theta^L\\
\begin{split}
\ddot \varphi^L 
&= -(
\ddot n^L_3 \cos \theta^L 
+ \left(\dot \varphi^L\right)^2 \cos \varphi^L \cos^2 \theta^L 
+ \left(\dot \theta^L\right)^2 \cos \varphi^L
- \ddot n^L_2 \cos \varphi^L \sin \theta^L \\
&- 2 \dot \varphi^L \dot \theta^L \cos \theta^L \sin \varphi^L \sin \theta^L
)  / \sin \varphi^L \cos^2 \theta^L. 
\end{split}
\end{cases}
\end{equation}
Again, each of these conversion relations suffers singularities at large rotations similar to that of 2D formulations.
To eliminate the singularities, 
we first apply the weighted average to the conversion between
between $\ddot \theta^L$ and $\bm{\ddot n}^L$ 
with Eqs. \eqref{3D_n_second_derivative1} and \eqref{3D_n_second_derivative2} as 
\begin{equation}\label{3D_conversion_relation1}
\begin{split}
\ddot \theta^L 
&= -\left(
\ddot n^L_3 \cos \varphi^L+\ddot n^L_1 \sin \varphi^L 
+ \left(\dot \varphi^L \right)^2 \cos \theta^L + \left(\dot \theta^L \right)^2 \cos \theta^L
\right) \sin \theta^L \\
&+\left(
\sin \theta^L \left(\dot \theta^L \right)^2 - \ddot n^L_2
\right)  \cos \theta^L.
\end{split}
\end{equation}
Then, for the conversion relation between $\ddot \varphi^L$ and $\bm{\ddot n}^L$, 
according to Eq. \eqref{3D_n_second_derivative1}, 
we can rewrite the relation as  
\begin{equation} \label{3D_n_second_derivative4}
\cos \theta^L = B / \ddot \varphi^L, 
\end{equation}
where $B$ denotes the numerator of 
the 2nd expression in Eq. \eqref{3D_n_second_derivative1}. 
We further denote the numerators of the 2nd expressions 
in Eqs. \eqref{3D_n_second_derivative2}
and \eqref{3D_n_second_derivative3}, respectively, as $B_1$ and $B_2$.
Inserting Eq. \eqref{3D_n_second_derivative4} into 
Eqs. \eqref{3D_n_second_derivative2} and \eqref{3D_n_second_derivative3}, 
we have 
\begin{equation}\label{3D_n_second_derivative5}
\begin{cases}
\ddot \varphi^L = \left(B^2 \cos \varphi^L\right) / B_1\\
\ddot \varphi^L = \left(B^2 \sin \varphi^L\right) / B_2,
\end{cases}
\end{equation}
and obtain the weighted average of the conversion relation as 
\begin{equation}\label{3D_conversion_relation2}
\ddot \varphi^L = \frac{B_1 B^2 \cos \varphi^L + B_2 B^2 \sin \varphi^L}
{B_1^2 + B_2^2}.
\end{equation}
Again, one can easily find that the present relations recover 
$\bm{\ddot{\theta}}^L = (\ddot{\theta}^L, \ddot{\varphi}^L)  = (-\ddot{n}_2^L, \ddot{n}_1^L)$ for small rotations.
\subsection{Time stepping} 
For the time integration of plate and shell dynamics, 
we use the position-based Verlet scheme \cite{zhang2021multi}. 
At the beginning of each time step, 
the deformation gradient, 
particle position, rotation angles and pseudo normal
are updated to the midpoint of the $n$-th time step as 
\begin{equation}\label{eq:verlet-first-half-solid}
\begin{cases}
\mathbb{F}^{L, n + \frac{1}{2}} = \mathbb{F}^{L, n} + \frac{1}{2} \Delta t \dot{\mathbb{F}}^{L, n}\\
\bm{r}_m^{n + \frac{1}{2}} = \bm{r}_m^n + \frac{1}{2} \Delta t \bm{\dot u}_m^n\\
\bm{\theta}^{L, n + \frac{1}{2}} = \bm{\theta}^{L, n} + \frac{1}{2} \Delta t \bm{\dot \theta}^{L, n}\\
\bm{n}^{L, n + \frac{1}{2}} = \bm{n}^{L, n} + \frac{1}{2} \Delta t \bm{\dot n}^{L, n}.
\end{cases}
\end{equation}
After the stress correction and Gauss-Legendre quadrature 
Eqs. \eqref{gaussian_quadrature_rule1} and \eqref{gaussian_quadrature_rule2}, 
the conservation equations are solved to 
obtain the $\bm{\ddot u}_m^{n+1}$ and $\bm{\ddot n}^{n+1}$. 
After transforming $\bm{\ddot n}^{n+1}$ to $\bm{\ddot n}^{L, n+1}$, 
$\bm{\ddot \theta}^{L, n+1}$ is obtained 
through the conversion relation between the pseudo normal and rotation angle, 
i.e., 
Eq. \eqref{2D_conversion_relation3} for 2D problems
and
Eqs. \eqref{3D_conversion_relation1} and \eqref{3D_conversion_relation2} 
for 3D problems.
At this point, 
the velocity and angular velocity are updated by 
\begin{equation}\label{eq:verlet-first-mediate-solid}
\begin{cases}
\bm{\dot u}_m^{n + 1} = \bm{\dot u}_m^{n} +  \Delta t  \bm{\ddot u}_m^{n+1}\\
\bm{\dot \theta}^{L, n + 1} = \bm{\dot \theta}^{L, n} +  \Delta t  \bm{\ddot \theta}^{L, n+1},
\end{cases}
\end{equation}
and the change rate of pseudo normal $\bm{\dot {n}} $ is updated 
by Eq. \eqref{2D_n_first_derivative} or \eqref{3D_n_first_derivative}.
Finally, 
the change rate of the deformation gradient $\dot{\mathbb{F}}^{L, n + 1}$ 
is estimated by Eq. \eqref{change_rate_deformation_gradient}, 
and then the deformation gradient, density, 
particle position, rotation angles and pseudo normal 
are updated to the new time step with 
\begin{equation}\label{eq:verlet-second-final-solid}
\begin{cases}
\mathbb{F}^{L, n + 1} = \mathbb{F}^{L, n + \frac{1}{2}} + \frac{1}{2} \Delta t \dot{\mathbb{F}}^{L, n + 1}\\
\rho^{n + 1} = \left(J_m^{n + 1} \right)^{-1} \rho^0  \\
\bm{r}_m^{n + 1} = \bm{r}_m^{n + \frac{1}{2}} + \frac{1}{2} \Delta t \bm{\dot u}_m^{n + 1}\\
\bm{\theta}^{L, n + 1} = \bm{\theta}^{L, n + \frac{1}{2}} + \frac{1}{2} \Delta t \bm{\dot \theta}^{L, n + 1}\\
\bm{n}^{L, n + 1} = \bm{\theta}^{L, n + \frac{1}{2}} + \frac{1}{2} \Delta t \bm{\dot n}^{L, n + 1}.
\end{cases}
\end{equation}
For the numerical stability, the time-step size $\Delta t$ is given by
\begin{equation}\label{eq:dt}
\Delta t   =  \text{CFL}\min\left(\Delta t_1, \Delta t_2, \Delta t_3 \right),
\end{equation}
where
\begin{equation}\label{time_step_size}
\begin{cases}
	\Delta t_1   =  \min\left(\frac{h}{c_v + |\bm{\dot u}_m|_{max}},
	\sqrt{\frac{h}{|{\bm{\ddot u}_m}|_{max}}} \right)\\
	\Delta t_2   =  \min\left(\frac{1}{c_v + |\bm{\dot \theta}^L|_{max}},
	\sqrt{\frac{1}{|{\bm{\ddot \theta}^L}|_{max}}} \right)\\
	\Delta t_3   =  h \left( \frac{\rho \left(1 - \nu^2 \right) / E}  {2 + \left(\pi^2/12 \right) \left(1 - \nu \right) \left[ 1 + 1.5 \left(h/d \right)^2  \right]   } 	\right)^{1/2}.\\
\end{cases}
\end{equation}
Note that the time-step size $\Delta t_3$ 
follows the Refs. \cite{lin2014efficient, tsui1971stability} 
and depends on the thickness and material properties, 
and the present Courant-Friedrichs-Lewy (CFL) number is set as 0.6.
An overview of the complete solution procedure 
is presented in Algorithm \ref{algorithm1}.
{
	\linespread{1.0} \selectfont
	\begin{algorithm}[htb!]
		\small{
		Setup parameters and initialize the simulation\;
		Construct the particle-neighbor list and compute the kernel values\;
		Compute the correction matrices $\widetilde{\mathbb{B}}^{0, \bm{r}}$ and $\widetilde{\mathbb{B}}^{0, \bm{n}}$ for each particle (Section \ref{sec:consistency correction})\;
		\While{simulation is not finished}
		{
			Compute the time-step size $\Delta t$ using Eq. \eqref{eq:dt}\;
			Update the deformation gradient tensor $\mathbb{F}^L$, particle position $\bm{r}_m$, rotation angle $\bm{\theta}^L$ and pseudo normal $\bm{n}$ for half time step $\Delta t/2$\; 
			Compute and correct the Cauchy stress $\fancy{$\sigma$}^l$ (Sections \ref{sec:constitutive_relation} and \ref{sec:stress_correction})\;
			Compute the resultants $\mathbb{N}^l$ and $\mathbb{M}^l$, and shear force $\bm{Q}^l$ (Eq. \eqref{resultants})\;
			Compute the acceleration $\bm{\ddot u}_m$ (Eqs. \eqref{discrete_dynamic_equation1} and \eqref{artificial-stress}) and $\bm{\ddot n}$ (Eqs. \eqref{discrete_dynamic_equation2} and \eqref{artificial-torque})\;
			Compute the angular acceleration $\bm{\ddot \theta}^L$ (Eq. \eqref{2D_conversion_relation3} for 2D problems, and Eqs. \eqref{3D_conversion_relation1} and \eqref{3D_conversion_relation2} for 3D problems)\;
			Update the velocity $\bm{\dot u}_m$ and angular velocity $\bm{\dot \theta}^L$ for a time step $\Delta t$\;
			Compute the change rate of pseudo normal $\bm{\dot n}^L$ using Eq. \eqref{2D_n_first_derivative} or \eqref{3D_n_first_derivative}\;
			Compute the change rate of the deformation gradient tensor $\partial \mathbb{F}^L / \partial t$ (Eq. \eqref{change_rate_deformation_gradient})\; 
			Update the deformation gradient tensor $\mathbb{F}^L$,  density $\rho$, particle position $\bm{r}_m$, 
			rotation angle $\bm{\theta}^L$ and pseudo normal $\bm{n}$ for another half time step $\Delta t/ 2$\;
		}
		Terminate the simulation.\
	}
		\caption{The present SPH method for plate/shell structures.}
		\label{algorithm1}
	\end{algorithm}
}

\section{Numerical examples}\label{sec:examples}
To demonstrate the accuracy and stability 
of the proposed surface-particle SPH method 
(denoted as shell method), 
this section investigates a series of benchmark tests 
where analytical or numerical reference data from the literature 
or/and volume-particle SPH method (denoted as volume method)
are available for qualitative and quantitative comparison. 
The smoothing length $h = 1.15~dp$, 
where $dp$ denotes the initial particle spacing, 
is employed in all the following simulations. 

\subsection{2D oscillating plate strip}
\begin{figure}[htb!]
	\centering
	\includegraphics[width = \textwidth]{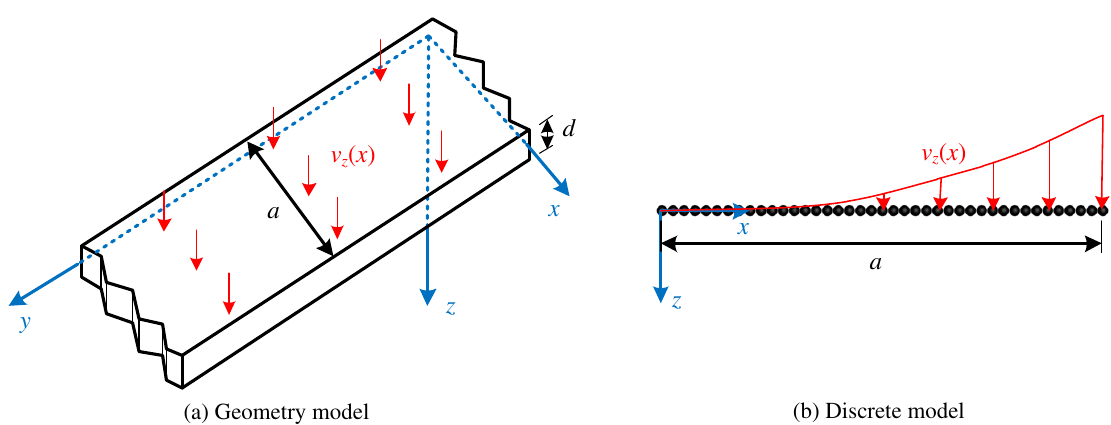}
	\caption{2D oscillating plate strip: Initial configuration with width $a = 0.2~\text{m}$.}
	\label{figs:2D_plate_setup}
\end{figure}

The first example involves a plate strip 
with initial uniform transverse  velocity 
along the length with one edge fixed and the others free, 
which has previously been theoretically \cite{landau1986course} 
and numerically \cite{gray2001sph, zhang2017generalized, wu2023essentially} 
investigated in the literature. 
As shown in Figure \ref{figs:2D_plate_setup}(a), 
this plate strip is assumed to be infinitely long along the $y$-axis 
with a finite width $a = 0.2~\text{m}$ along the $x$-axis. 
To demonstrate that 
both thin and moderately thick structures can be simulated, 
this plate strip is modeled with the thicknesses $d = 0.01~\text{m}$ and $0.001~\text{m}$.
The material properties are set as follows: 
density $\rho_0 = 1000.0 ~\text{kg} / \text{m}^3$, 
Young’s modulus $E = 2.0~\text{MPa}$, 
and Poisson’s ratio $\nu$ varies for different cases. 
Figure \ref{figs:2D_plate_setup}(b) shows 
the discrete model of the chosen cross-section with clamped edges at $x = 0$. 
The transverse velocity $v_z$ is applied to the plate strip as
\begin{equation}
	v_z(x) = v_f c \frac{f(x)}{f(a)},
\end{equation}
where $v_f$ is a constant that varies with different cases, 
and 
\begin{equation}
	\begin{split}
		f(x) &= \left(\sin(ka) + \sinh(ka) \right) \left(\cos(kx) - \cosh(kx) \right) \\
		& - \left(\cos(ka) + \cosh(ka) \right) \left(\sin(kx) - \sinh(kx) \right)
	\end{split}
\end{equation}
with $k$ determined by 
\begin{equation}
	\cos(ka) \cosh(ka) = -1
\end{equation}
and $ka = 1.875$. 
The frequency $\omega$ of the oscillating plate strip is theoretically given by 
\begin{equation}
	\omega ^2 = \frac{E d^2 k^4}{12 \rho \left(1 - \nu^2 \right)}.
\end{equation}

\begin{figure}[htb!]
	\centering
	\includegraphics[trim = 2mm 6mm 2mm 2mm, width=\textwidth] {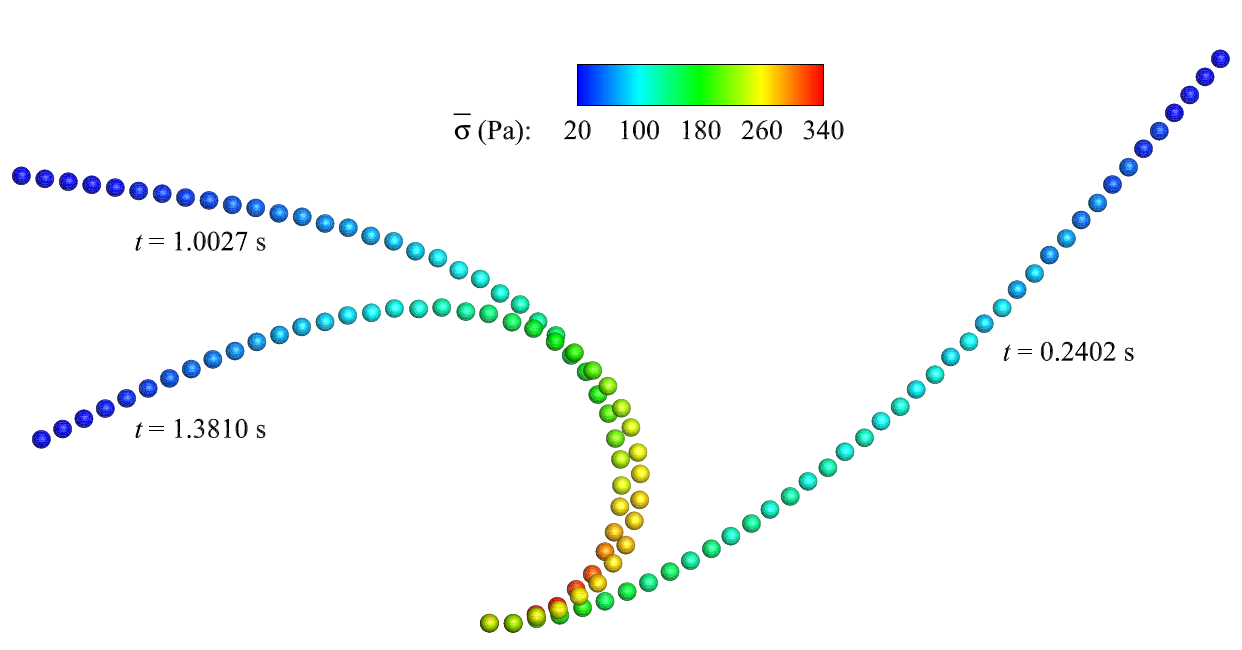}
	\caption{2D oscillating plate strip: Deformed particle configuration colored by von Mises stress $\bar\sigma$ of the mid-surface at serial time instants
	with the width $a = 0.2~\text{m}$, thickness $d = 0.02~\text{m}$, $v_f = 0.01$, and spatial particle resolution $a / dp = 40$. 
		The material is modeled with density $\rho_0 = 1000.0 ~\text{kg} / \text{m}^3$, Young’s modulus $E = 2.0~\text{MPa}$, and Poisson’s ratio $\nu = 0.4$. }
	\label{figs:2D_oscillating_plate_stress}
\end{figure}

\begin{figure}[htb!]
	\centering
	\includegraphics[trim = 2mm 6mm 2mm 2mm, width=\textwidth] {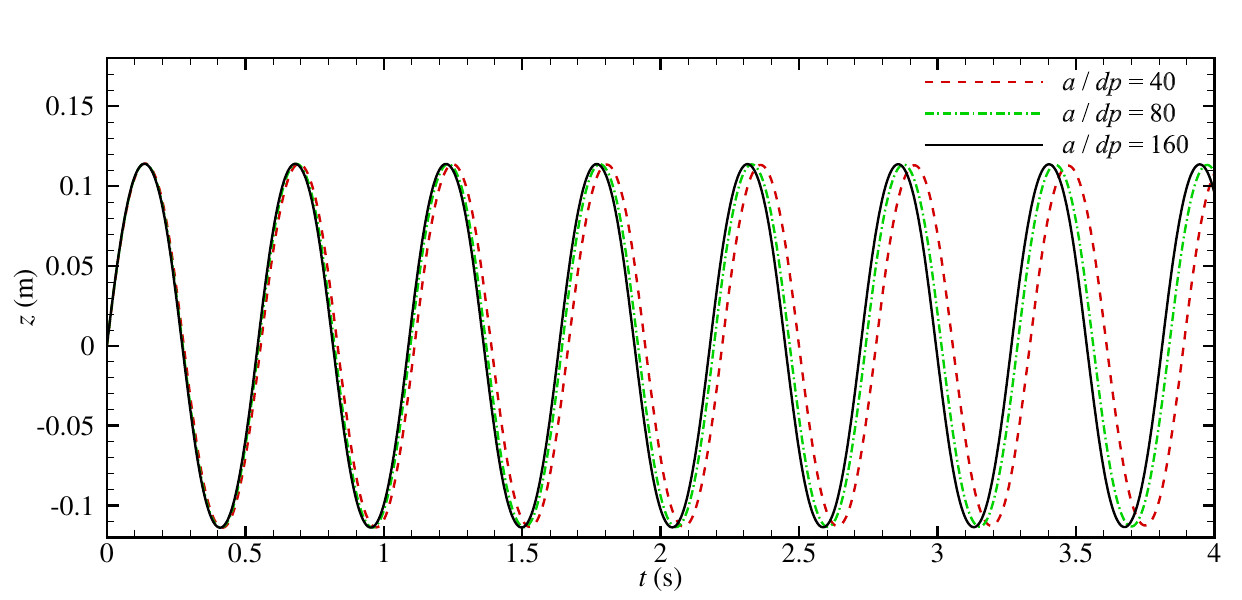}
	\caption{2D oscillating plate strip: Time history of the vertical position $z$ 
		observed at the plate strip endpoint when $d = 0.01~\text{m}$, $v_f = 0.025$ and $\nu = 0.4$.}
	\label{figs:2D_oscillating_plate_convergence}
\end{figure}

Figure \ref{figs:2D_oscillating_plate_stress} shows the particles 
with von Mises stress $\bar\sigma$ contour for the case of 
$d = 0.001~\text{m}$, $v_f = 0.01$, $\nu = 0.4$,  
and the initial particle spacing $dp = a / 40 = 0.005~\text{m}$.
It should be noted that 
the present method predicts smooth deformation and stress fields 
without singularities for large rotations (more than $\pi$). 
Three different spatial resolutions,
$a / dp =40$, $a / dp =80$, and $a / dp =160$,
are tested in the convergence study. 
Figure \ref{figs:2D_oscillating_plate_convergence} shows 
the time history of vertical position $z$ of the strip endpoint
with $d = 0.01~\text{m}$, $v_f = 0.025$ and $\nu = 0.4$. 
It can be observed that typical 2nd-order convergence has been achieved. 
In addition, 
a long-term simulation is performed herein to demonstrate the numerical stability of the proposed formulation.  
For quantitative validation, 
\begin{table}[tb!]
	\centering
	\caption{2D oscillating plate strip: Quantitative validation of the oscillation period for $a = 0.2~\text{m}$ and $d = 0.01~\text{m}$ with various $v_f $ and $\nu$.}
	\begin{tabular}{ccccc}
		\hline
		$v_f $   & $\nu$  & $T_\text{Shell model}$ & $T_\text{Theoretical}$ & Error \\ 
		\hline
		0.025		   & 0.22  	   & 0.58137     & 0.54018      & 7.63\%\\
		0.05 			& 0.22  	& 0.57715     & 0.54018      & 6.92\%\\
		0.1	 		     & 0.22      & 0.56801     & 0.54018      & 5.15\%\\
		~\\
		0.025		   & 0.30  	   & 0.56804     & 0.52824      & 7.53\%\\
		0.05 			& 0.30  	& 0.56308     & 0.52824      & 6.60\%\\
		0.1	 		     & 0.30      & 0.55481     & 0.52824      & 5.03\%\\
		~\\
		0.025		   & 0.40  	   & 0.54447     & 0.50752      & 7.28\%\\
		0.05 			& 0.40  	& 0.53683     & 0.50752      & 5.78\%\\
		0.1	 		     & 0.40      & 0.53252     & 0.50752      & 4.93\%\\
		\hline	
	\end{tabular}
	\label{tab:oscillating_plate_period1}
\end{table}
\begin{table}[tb!]
	\centering
	\caption{2D oscillating plate strip: Quantitative validation of the oscillation period for $a = 0.2~\text{m}$ and $d = 0.001~\text{m}$ with various $v_f $ and $\nu$.}
	\begin{tabular}{ccccc}
		\hline
		$v_f $   & $\nu$  & $T_\text{Shell model}$ & $T_\text{Theoretical}$ & Error \\ 
		\hline
		0.0025		   & 0.22  	   & 5.80249     & 5.40182      & 7.42\%\\
		0.005 			& 0.22  	& 5.75544     & 5.40182      & 6.55\%\\
		0.01	 		 & 0.22      & 5.64181     & 5.40182      & 4.44\%\\
		~\\
		0.0025		   & 0.30  	   & 5.66756     & 5.28243      & 7.29\%\\
		0.005 			& 0.30  	& 5.61006     & 5.28243      & 6.20\%\\
		0.01	 		 & 0.30      & 5.49156     & 5.28243      & 3.96\%\\
		~\\
		0.0025		   & 0.40  	  & 5.42826     & 5.07519      & 6.96\%\\
		0.005 			& 0.40     & 5.34224     & 5.07519      & 5.26\%\\
		0.01	 		 & 0.40     & 5.27522     & 5.07519      & 3.94\%\\
		\hline	
	\end{tabular}
	\label{tab:oscillating_plate_period2}
\end{table}
Tables \ref{tab:oscillating_plate_period1} and \ref{tab:oscillating_plate_period2} 
detail the oscillation period $T$ for a wide range of $v_f$ and $\nu$, 
obtained by the present method with the spatial particle resolution $a / dp =160$, 
when thickness $d = 0.01~\text{m}$ and $0.001~\text{m}$,
respectively, 
and the comparison to theoretical solution obtained form small perturbation analysis. 
The differences, 
which are less than 8.00\% for $\nu = 0.22$ and decrease to about 5.00\% when the Poisson’s ratio is increased to 0.4, 
are acceptable.
\subsection{3D square plate}
\label{3D_plate_steady}
\begin{figure}[tb!]
	\centering
	\includegraphics[trim = 2mm 6mm 2mm 2mm, width =0.85 \textwidth]{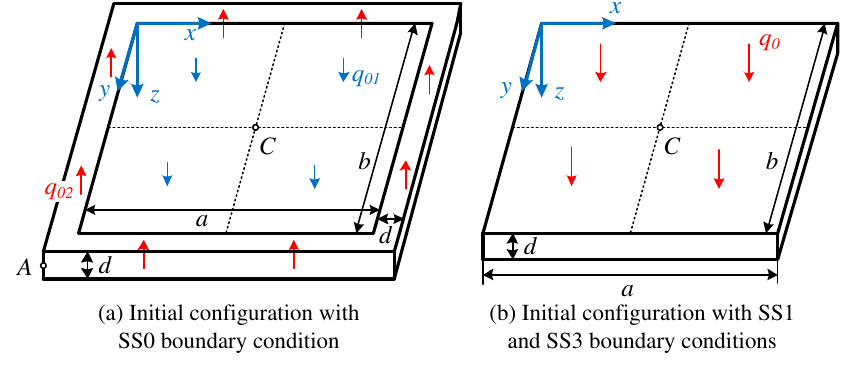}
	\caption{3D square plate: Problem setup with $a = b = 254~\text{mm}$ and thickness $d = 25.4~\text{mm}$.}
	\label{figs:3D_square_plate_setup}
\end{figure}

In this section, 
a 3D square plate under different types of boundary conditions 
is considered for quasi-steady analyses, 
as shown in Figure \ref{figs:3D_square_plate_setup}. 
With side length $a = b = $ 254~mm 
and thickness $d =$ 25.4~mm, 
the plate material is defined with 
density $\rho_0 = 1600 ~\text{kg} / \text{m}^3$, 
Young’s modulus $E = 53.7791 ~\text{GPa}$ 
and Poisson’s ratio $\nu = 0.3$. 
Three types of boundary conditions denoted as 
SS0, SS1 and SS3 following 
Refs. \cite{reddy2006theory, lin2014efficient} are implemented as
\begin{itemize}
     	\item[$\blacktriangleright$] SS0: \qquad constrained mass center without constrained boundaries;
		\item[$\blacktriangleright$] SS1: \qquad  $u = w = \varphi = 0$ on edges parallel to $x$-axis 
		and $v = w = \theta = 0$ on edges parallel to $y$-axis;
	    \item[$\blacktriangleright$] SS3: \qquad $u = v = w = 0$ on all edges.
\end{itemize}
Note that,
for the case of SS0,
the outer square ring with width $d$ 
is imposed with negative pressure $q_{02}$. 
The uniformly distributed loads are parameterized by 
the loading factors $\bar P$ and $\bar P_1$ 
as $q_0  = \bar PE\left( {d/a} \right)^4$,  
$q_{01}  = \bar P_1 E\left( {d/a} \right)^4$
and $q_{02} (2ad + 2bd + 4 d^2) = q_{01} a b$, 
so that the applied negative force along the $z$-axis 
prevents the center of mass from moving. 

\begin{figure}[htb!]
	\centering
	\includegraphics[trim = 2mm 6mm 2mm 2mm, width=\textwidth] {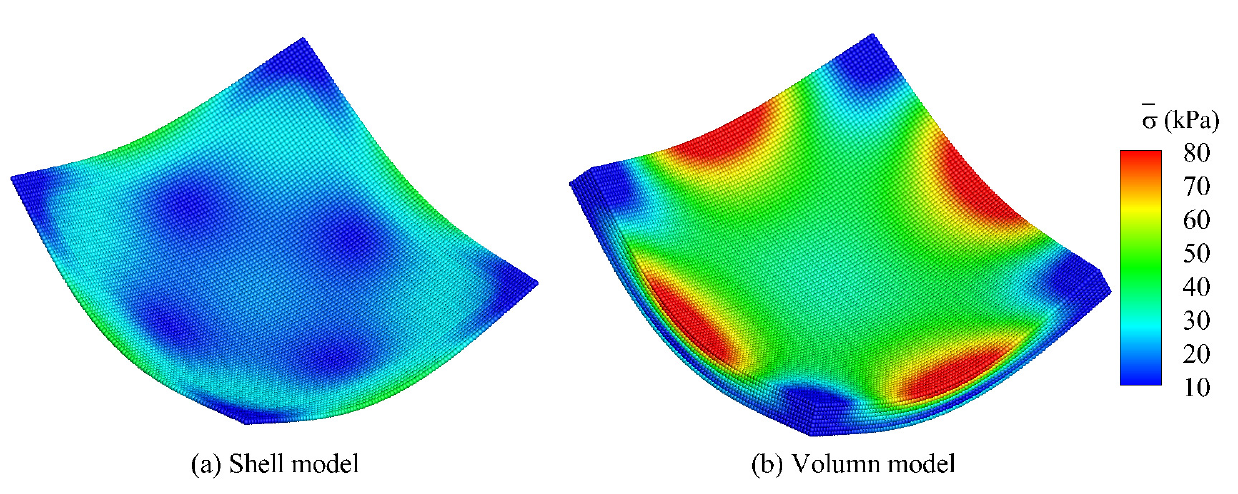}
	\caption{3D square plate: Particles
		colored by von Mises stress $\bar\sigma$ 
		of tests with SS0 obtained by the present shell (left) 
		and volume (right) methods
		under the loading factor $\bar P_1 = 25$. 
		Note that the left panel shows the stress $\bar\sigma$ of the plate mid-surface.
		The material is modeled with the density $\rho_0 = 1600 ~\text{kg} / \text{m}^3$,  
		Young’s modulus $E = 53.7791 ~\text{GPa}$ and Poisson’s ratio $\nu = 0.3$.
		The spatial particle resolution is set as $d / dp = 8$.}
	\label{figs:3D_square_plate_comparison_volumn_stress}
\end{figure}

For comprehensive validation, 
a convergence study of tests with SS0 is conducted, 
and the results are compared with those obtained by the volume method 
released in the SPHinXsys repository \cite{zhang2021sphinxsys}. 
Figure \ref{figs:3D_square_plate_comparison_volumn_stress} 
shows the particle distribution and stress fileds 
under the loading factor $\bar P_1 = 25$ 
with the spatial discretization $d/dp = 8$. 
\begin{figure}[htb!]
	\centering
	\includegraphics[trim = 2mm 6mm 2mm 2mm, width=\textwidth] {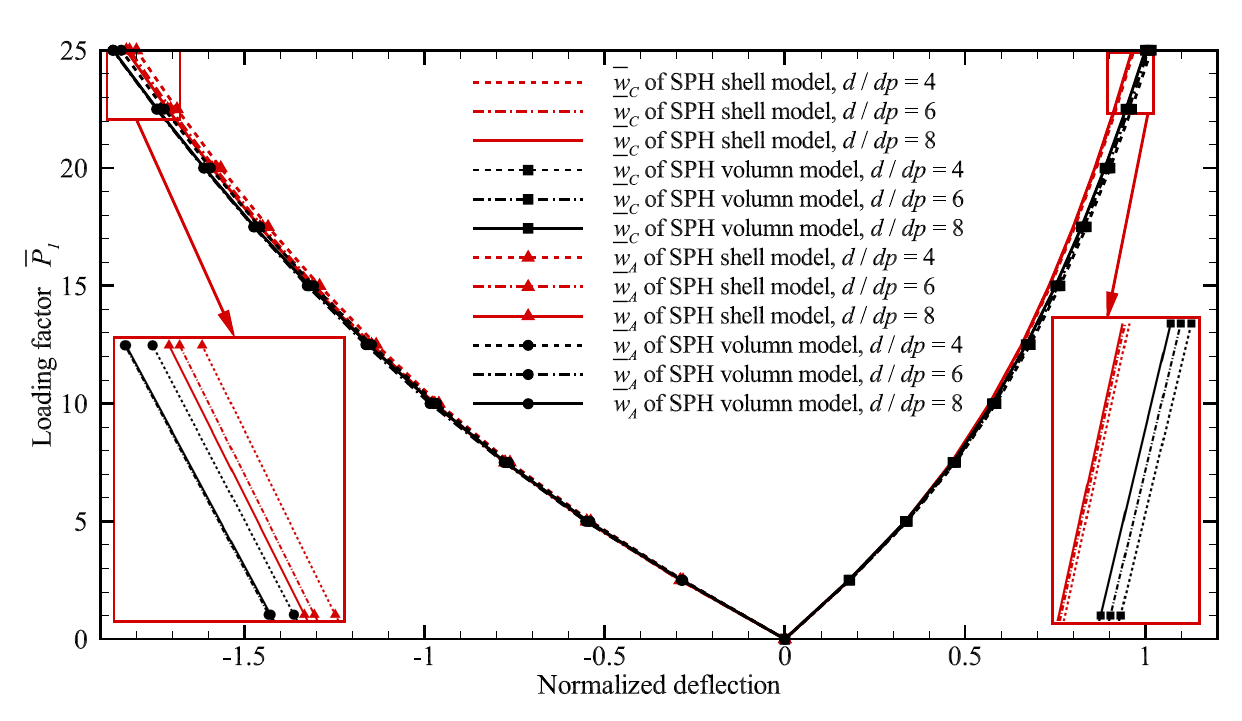}
	\caption{3D square plate: Load-deflection curves of tests with SS0 
		under three different spatial resolutions, 
		and their comparison with those of the volume method \cite{zhang2021sphinxsys}.}
	\label{figs:3D_square_plate_comparison_volumn}
\end{figure}
Figure \ref{figs:3D_square_plate_comparison_volumn} shows 
the non-dimensional deflection $\bar w_C  = w_C /d$ 
and $\bar w_A = w_A / d$  
probed at the central point $C$ 
and corner point $A$, respectively, 
obtained by both SPH shell and volume methods.
It should be emphasized that there are only quite small differences 
between the results of the present reduced-dimensional and full-dimensional models. 

\begin{figure}[htb!]
	\centering
	\subfigure[SS1 boundary conditon]{
		\includegraphics[trim = 2mm 6mm 2mm 4mm,  width=\linewidth]{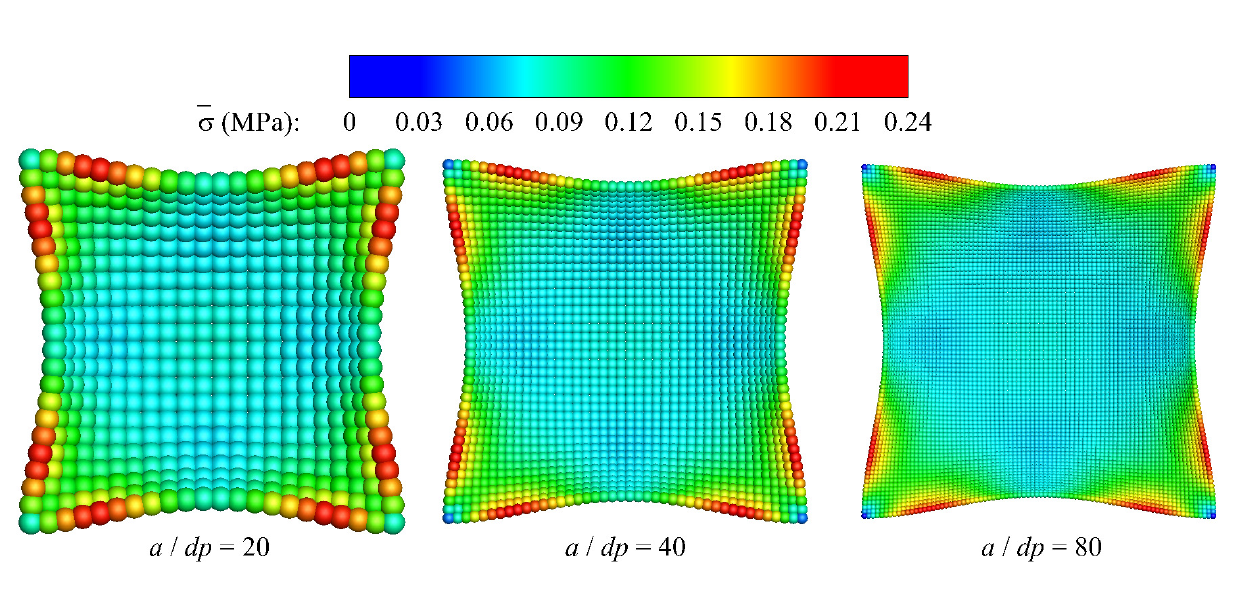}}
	\subfigure[SS3 boundary conditon]{
		\includegraphics[trim = 2mm 6mm 2mm 4mm,  width=\linewidth]{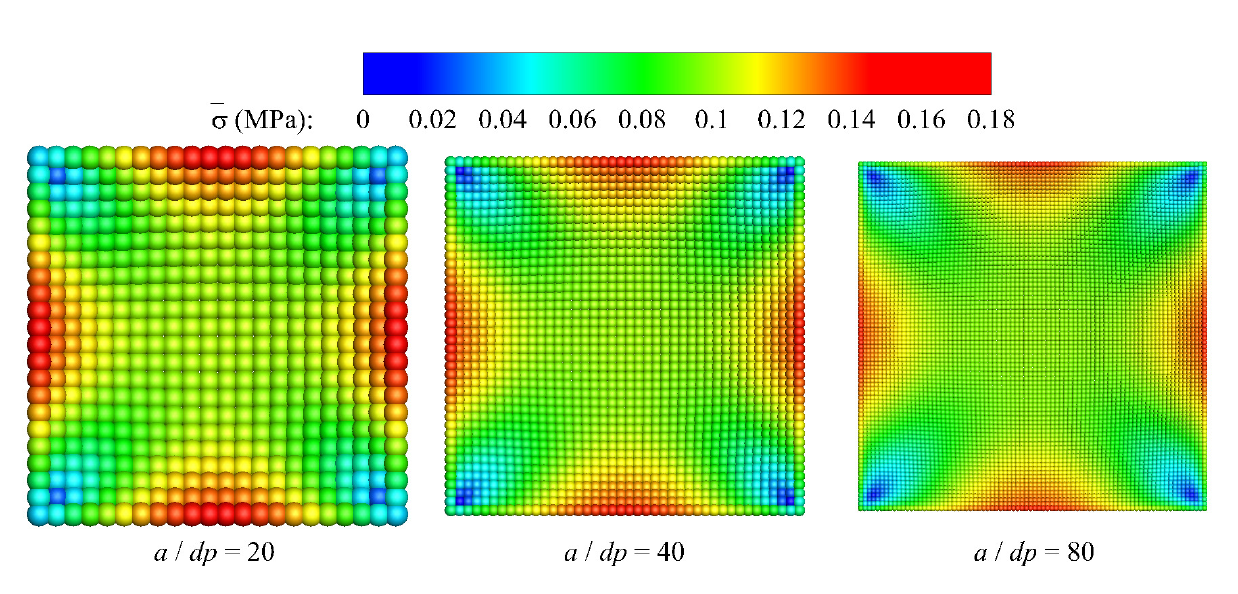}}
	\caption{3D square plate: Particles
		colored by von Mises stress $\bar\sigma$ of the mid-surface 
		with particle refinement under the loading factor $\bar P = 200$ 
		and SS1 and SS3 boundary conditions. 
		The plate material has parameters of the density $\rho_0 = 1600 ~\text{kg} / \text{m}^3$,  
		Young’s modulus $E = 53.7791 ~\text{GPa}$ and Poisson’s ratio $\nu = 0.3$. }
	\label{figs:3D_square_plate_SS1_SS3_stress}
\end{figure}

The particles colored by von Mises stress $\bar\sigma$ at the mid-surface
for three spatial discretizations, 
$a/dp = 20$, $a/dp = 40$ and $a/dp = 80$, 
with the SS1 and SS3 boundary conditions 
under $\bar P = 200$ 
are shown in Figure \ref{figs:3D_square_plate_SS1_SS3_stress}. 
\begin{figure}[htb!]
	\centering
	\subfigure[SS1 boundary conditon]{
		\includegraphics[trim = 2mm 6mm 2mm 4mm,  width=0.48\linewidth]{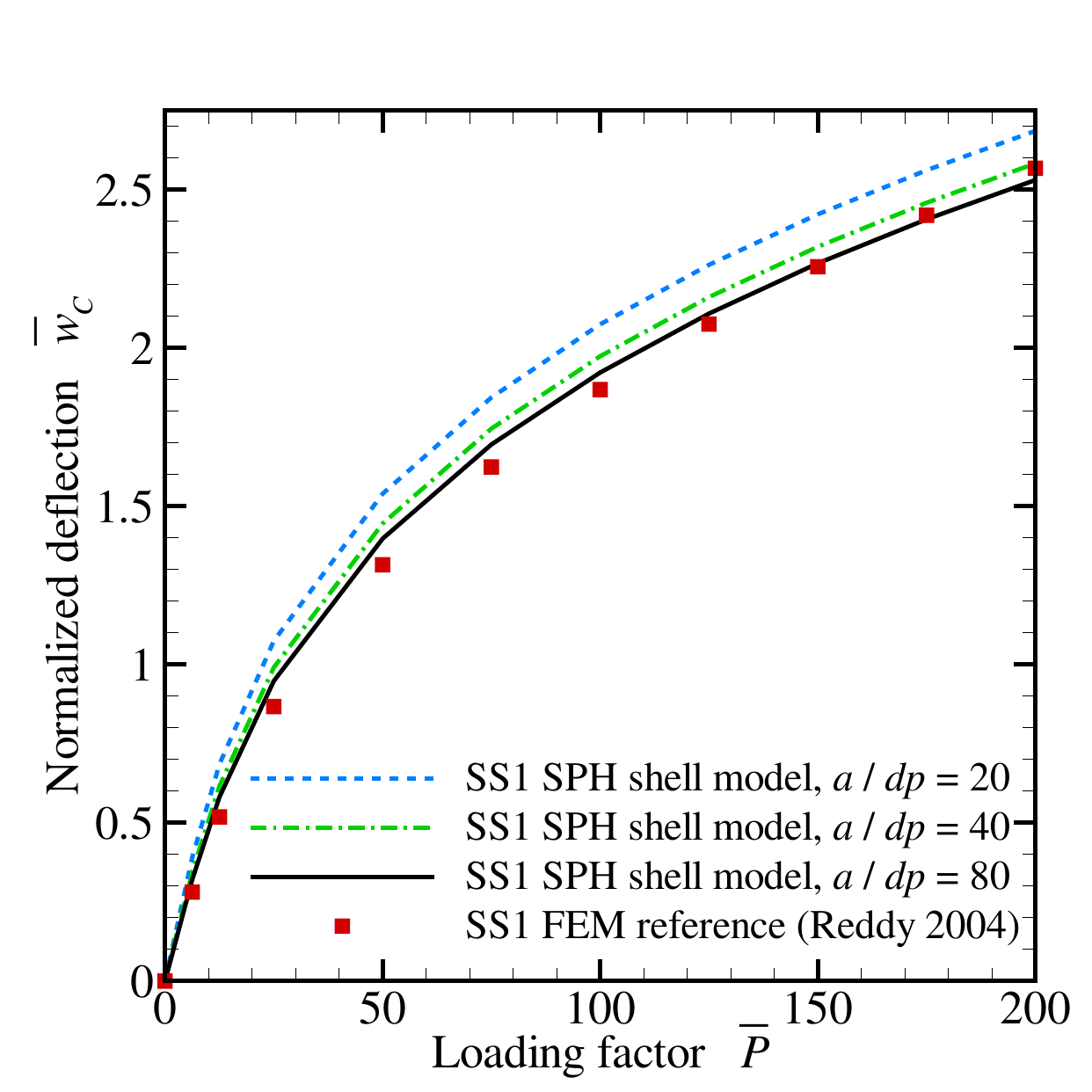}}
	\subfigure[SS3 boundary conditon]{
		\includegraphics[trim = 2mm 6mm 2mm 4mm,  width=0.48\linewidth]{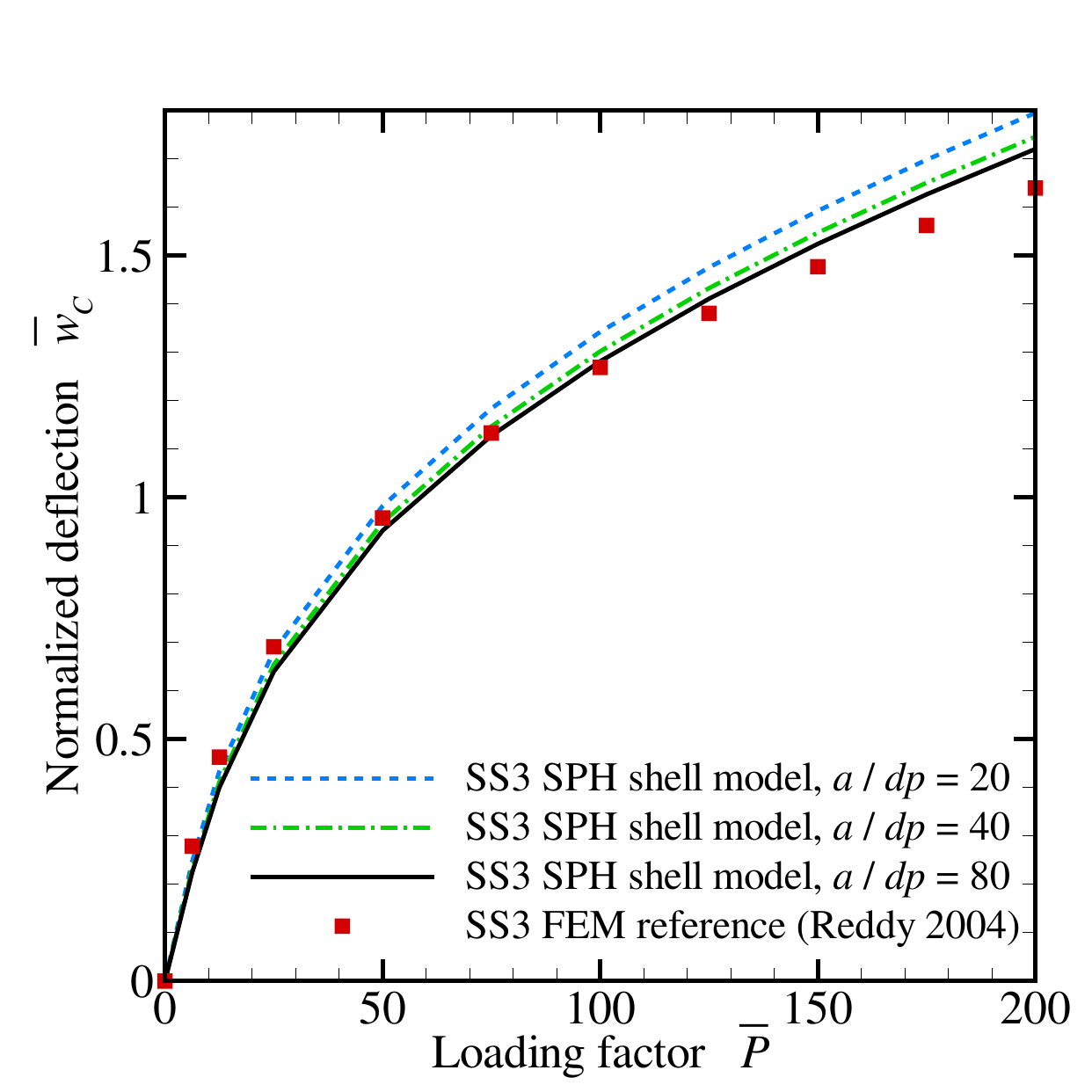}}
	\caption{3D square plate: Load-deflection curves of tests with SS1 and SS3 
		under three different spatial resolutions, 
		and their comparison with that of Reddy \cite{reddy2006theory}.}
	\label{figs:3D_square_plate_SS1_SS3_comparison}
\end{figure}
It can be observed that 
the regular particle distribution and smooth stress field are obtained. 
Also, both the deformation and von Mises stress $\bar \sigma $ 
exhibit good convergence properties
with particle refinement. 
In order to demonstrate the accuracy of the present method, 
the non-dimensional deflections $\bar w_C$ for tests with
SS1 and SS3 under various spatial resolutions 
are compared to those of the Ref. \cite{reddy2006theory}.
As shown in Figs. \ref{figs:3D_square_plate_SS1_SS3_comparison}, 
the numerical results quickly converge to the reference solutions 
obtained by the finite element method (FEM) 
with increasing resolution.
\subsection{Dynamic response of a 3D square plate}
Following Ref. \cite{momenan2018new}, 
the 3D square plate studied in Section \ref{3D_plate_steady} 
is considered herein 
with the thickness $d = 12.7~\text{mm}$ 
and Young’s modulus $E = 68.94 ~\text{GPa}$.
The SS0 and SS3 boundary conditions are applied 
for dynamic analyses under a step loading of uniform normal pressure 
$q_{01} = q_0 = 2.068427 ~\text{MPa}$. 
For convergence study, 
three different spatial discretizations, 
i.e., $d/dp = 2$, $d/dp = 4$ and $d/dp = 8$, 
are considered.

\begin{figure}[htb!]
	\centering
	\includegraphics[trim = 2mm 6mm 2mm 2mm, width=\textwidth] {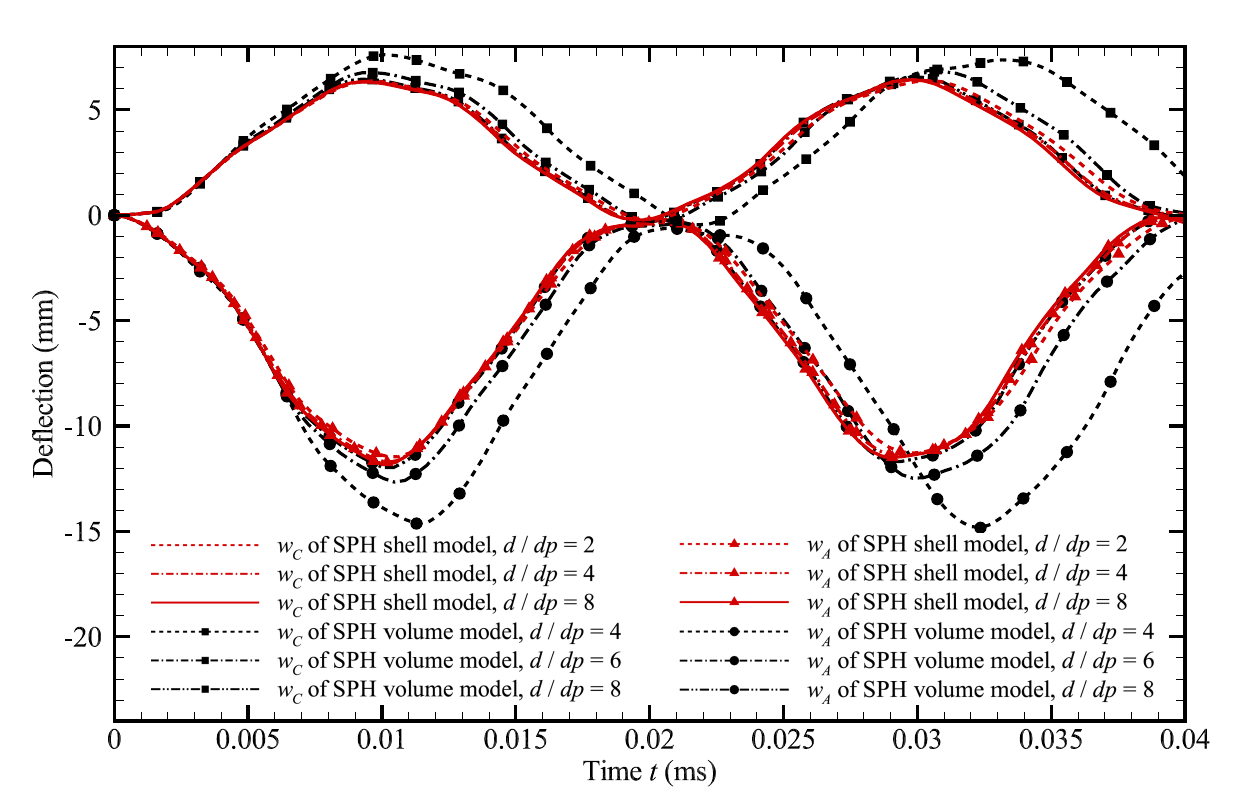}
	\caption{3D square plate with dynamic response: 
		Time history of the deflection $w_C$ and $w_A$ 
		probed at the central point $C$ and corner point $A$, respectively,  
		with SS0 boundary condition.
		The material is modeled with the density $\rho_0 = 1600 ~\text{kg} / \text{m}^3$,  
		Young’s modulus $E = 68.94~\text{GPa}$ and Poisson’s ratio $\nu = 0.3$.}
	\label{figs:3D_plate_dynamic_solution_comparison_volumn}
\end{figure}

\begin{figure}[htb!]
	\centering
	\includegraphics[trim = 2mm 6mm 2mm 2mm, width=0.5 \textwidth] {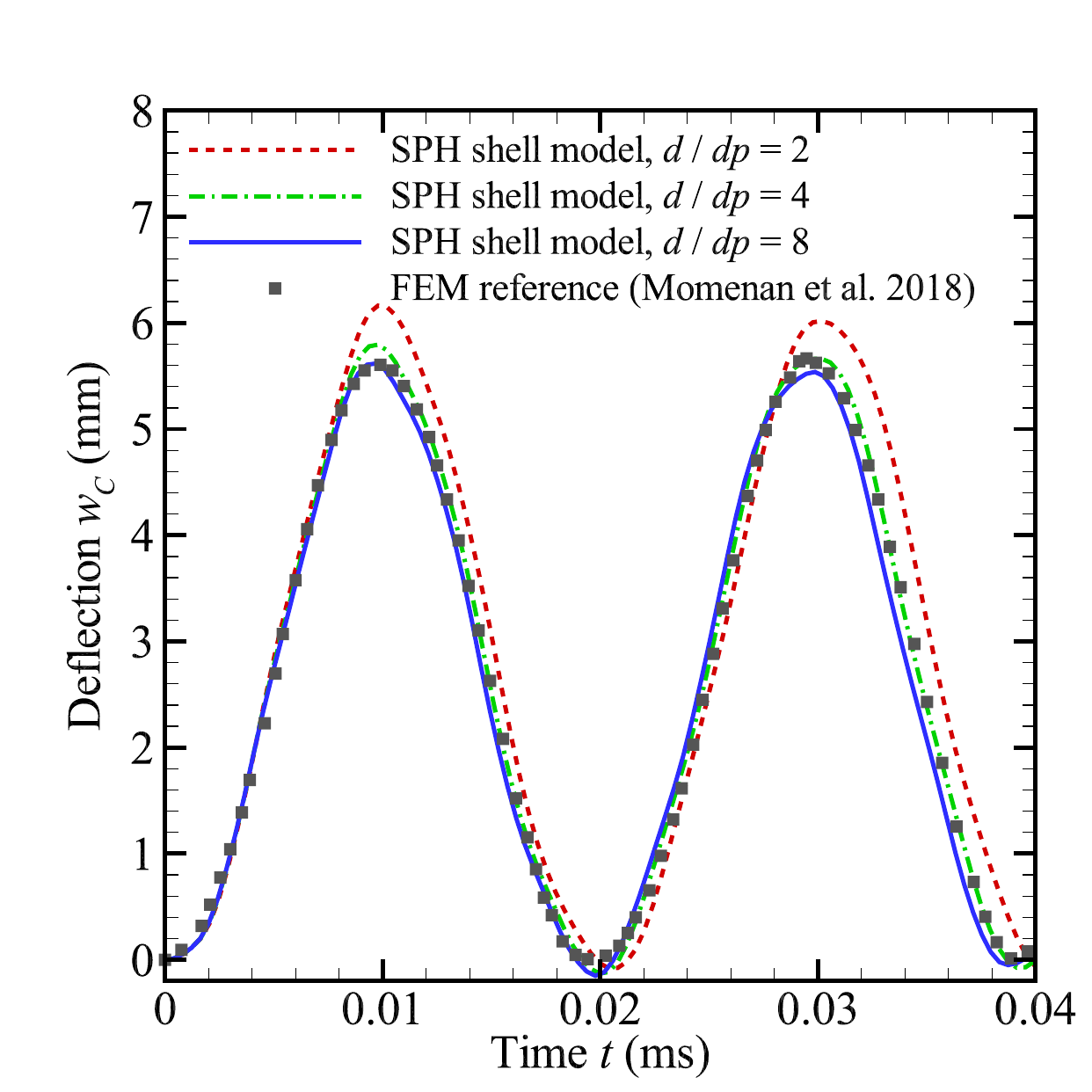}
	\caption{3D square plate with dynamic response: Time history of the deflection $w_C$ 
		observed at the central point $C$  
		with SS3 boundary condition.}
	\label{figs:3D_plate_dynamic_solution_comparison}
\end{figure}

For quantitative validation, 
Figure \ref{figs:3D_plate_dynamic_solution_comparison_volumn}
shows the time history of the deflections $w_C$ probed at the central point $C$ 
and $w_A$ at the corner point $A$
with SS0 boundary condition
and its comparison to the results obtained by the volume method. 
Also, Figure \ref{figs:3D_plate_dynamic_solution_comparison} shows
the time history of the deflection $w_C$ with the SS3 boundary condition
and its comparison with that of Ref. \cite{momenan2018new}. 
In general, 
the present results are in good agreements 
with those obtained by the volume method and of Ref. \cite{momenan2018new}.

\subsection{3D cantilevered plate}
Following Refs \cite{sze2002stabilized, sze2004popular, payette2014seven}, 
the static response of a 3D cantilevered plate 
subjected to a distributed end shear load $f_0$ is considered.
\begin{figure}[htb!]
	\centering
	\includegraphics[width = 0.5 \textwidth]{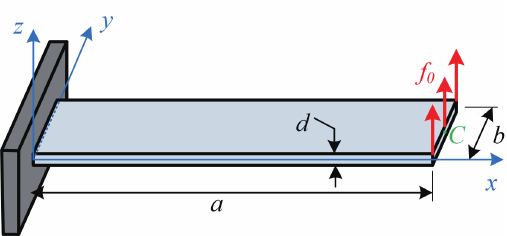}
	\caption{3D cantilevered plate: Initial configuration with the length $a = 10~\text{m}$, 
		width $b = 1~\text{m}$ and thickness $d = 0.1~\text{m}$.}
	\label{figs:3D_cantilevered_plate_setup}
\end{figure}
As shown in Figure \ref{figs:3D_cantilevered_plate_setup}, 
the plate with length $a = 10~\text{m}$, 
width $b = 1~\text{m}$ and thickness $d = 0.1~\text{m}$ is clamped at $y = 0$, 
and has material parameters of 
density $\rho_0 = 1100 ~\text{kg} / \text{m}^3$, 
Young’s modulus $E = 1.2 ~\text{MPa}$ and 
Poisson’s ratio $\nu = 0.0$. 
The shear load is parameterized by a loading factor 
$\bar F$ as $f_0  = \bar F EI / a^2$ 
with the inertia moment $I = bd^3/12$. 
Three different resolutions, 
i.e., $b/dp = 5$, $b/dp = 7$ and $b/dp = 9$, 
are considered for convergence study.

\begin{figure}[htb!]
	\centering
	\includegraphics[trim = 2mm 6mm 2mm 2mm, width=0.5\textwidth] {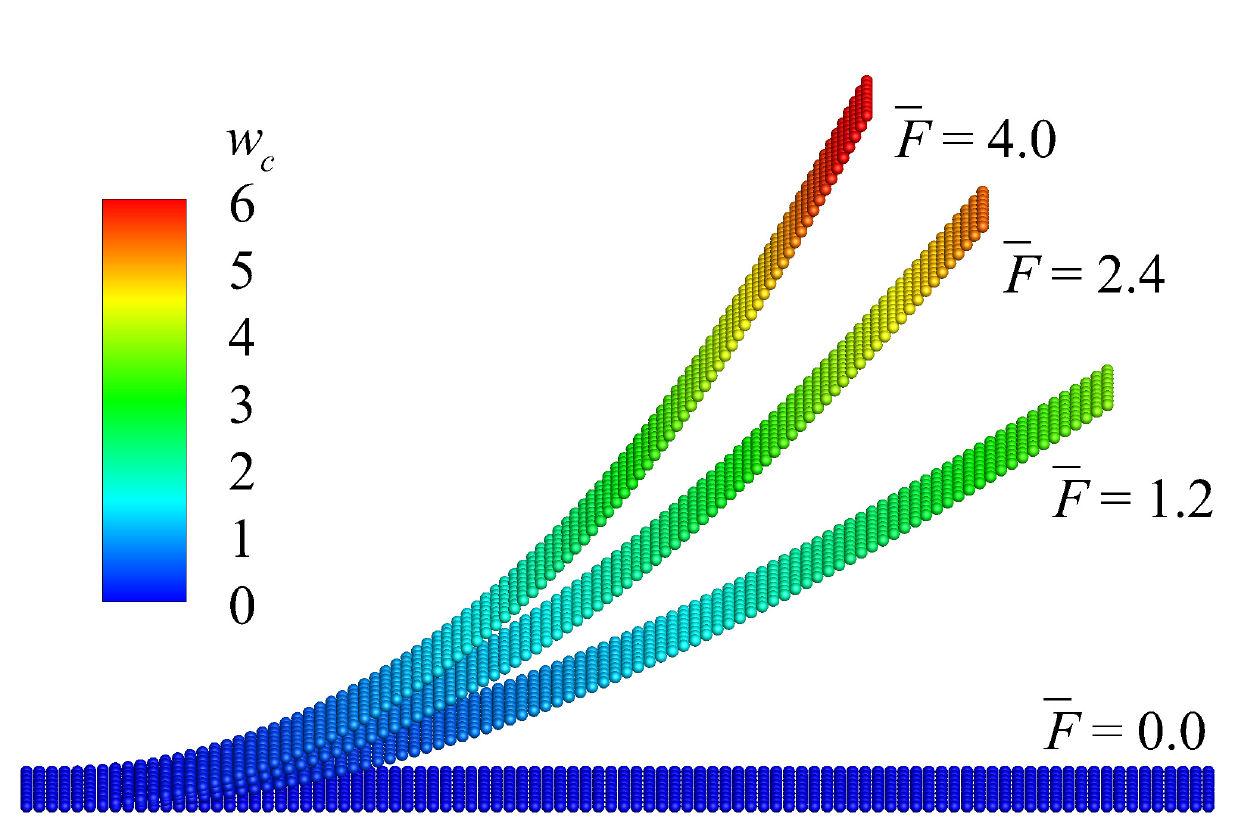}
	\caption{3D cantilevered plate: Particles 
		colored by the vertical displacement $w_c$ 
		under the various loading factor $\bar F$ at spatial resolution $b / dp = 9$.
		The material is set as the density $\rho_0 = 1100 ~\text{kg} / \text{m}^3$, 
		Young’s modulus $E = 1.2 ~\text{MPa}$ and Poisson’s ratio $\nu = 0.0$.}
	\label{figs:3D_cantilevered_plate_stress}
\end{figure}

Figure \ref{figs:3D_cantilevered_plate_stress} 
shows the particles 
colored by the vertical displacement
under different loading factor $\bar F$
at the spatial resolution $b / dp = 9$. 
A regular particle distribution and smooth vertical displacement field are noted.
\begin{figure}[htb!]
	\centering
	\includegraphics[trim = 2mm 6mm 2mm 2mm, width=0.5 \textwidth] {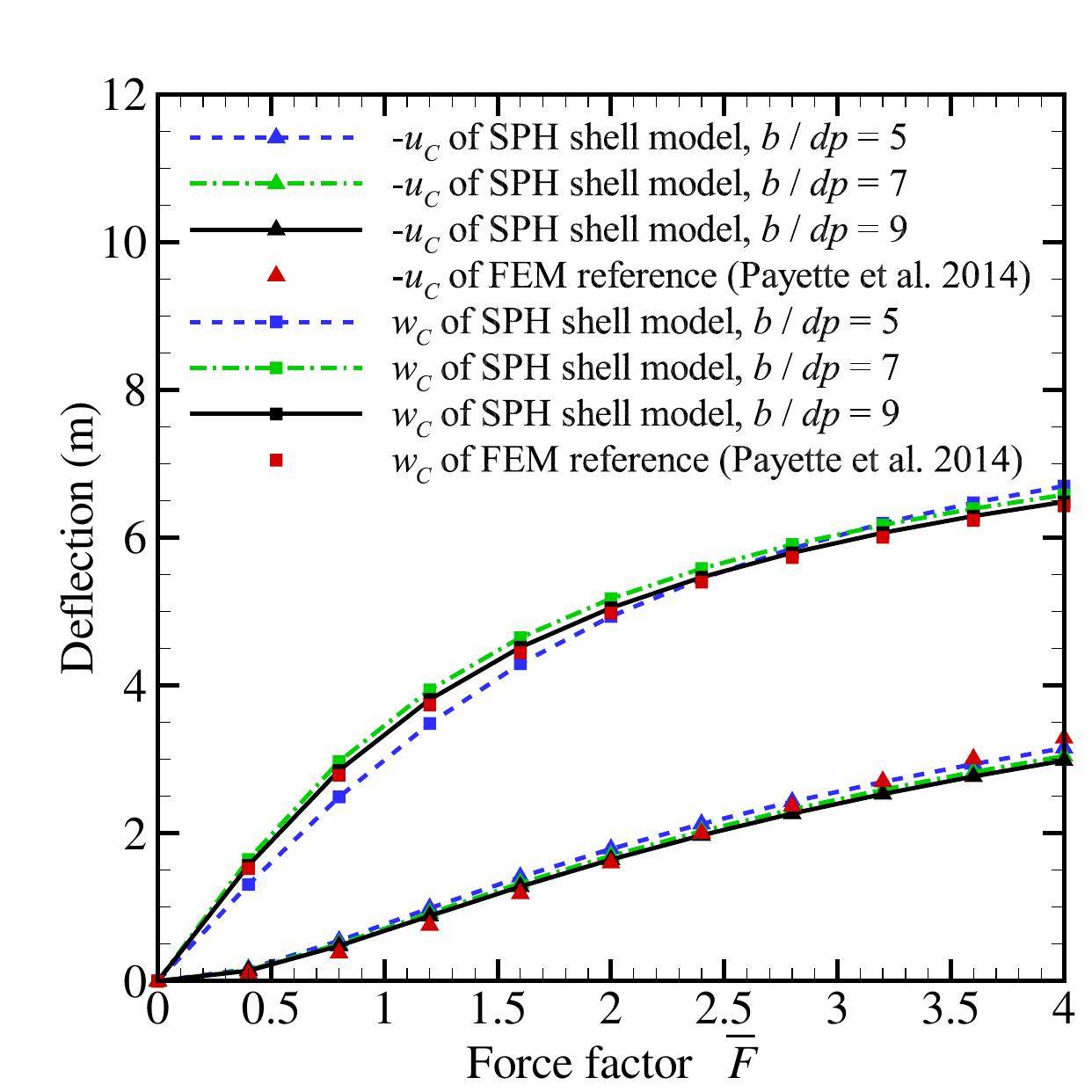}
	\caption{ 3D cantilevered plate: Load-deflection curves 
		with three various spatial discretizations, 
		and their comparison with that of Payette et al. \cite{payette2014seven}.}
	\label{figs:3D_cantilevered_plate_comparison}
\end{figure}
Figure \ref{figs:3D_cantilevered_plate_comparison} gives
the displacement $u_c$ and $w_c$ of the point $C$, 
defined in Figure \ref{figs:3D_cantilevered_plate_setup},  
as a function of the loading factor $\bar F$ and the initial particle spacing $dp$, 
and their comparison with those in Ref. \cite{payette2014seven}.
It can be noted that 
the displacement is converging rapidly,
again at about 2nd-order, 
with increasing resolution, 
demonstrating the accuracy of the present method.
\subsection{Scordelis-Lo roof}
\begin{figure}[htb!]
	\centering
	\includegraphics[width = 0.5 \textwidth]{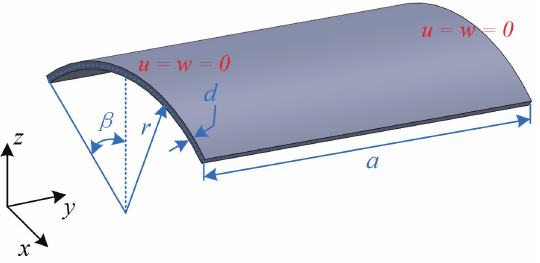}
	\caption{Scordelis-Lo roof: Initial configuration 
		with the length $a = 50~\text{m}$, 
		radius $r = 25~\text{m}$, 
		thickness $d = 0.25~\text{m}$ and $\beta = 40^\circ$.}
	\label{figs:3D_roof_setup}
\end{figure}

As shown in Figure \ref{figs:3D_roof_setup}, 
the Scordelis-Lo roof with length $a = 50~\text{m}$, 
radius $r = 25~\text{m}$, 
thickness $d = 0.25~\text{m}$ and $\beta = 40^\circ$
is considered herein, 
and the material properties are
density $\rho_0 = 36 ~\text{kg} / \text{m}^3$, 
Young’s modulus $E = 432 ~\text{MPa}$ and zero Poisson’s ratio. 
The roof is supported at its ends by fixed diaphragms, 
i.e. the translations in $x$ and $z$ directions are constrained, 
and subjected to a gravity loading of $g = 10 ~\text{m} / \text{s}^2$.

The FEM solution of the vertical displacement  $w$ 
at the midpoint of the side edge converges to 0.3024 $\text{m}$ 
as reported in Refs. \cite{belytschko1985stress, simo1989stress}. 
A sequentially refined resolutions 
of $b/dp = 15, 20, 25, 30~\text{and}~40$ 
with $b = 2 r \beta$ denoting the arc length of the roof end 
are considered to assess the convergence property 
of the present method. 
\begin{figure}[htb!]
	\centering
	\includegraphics[trim = 2mm 6mm 2mm 2mm, width=\textwidth] {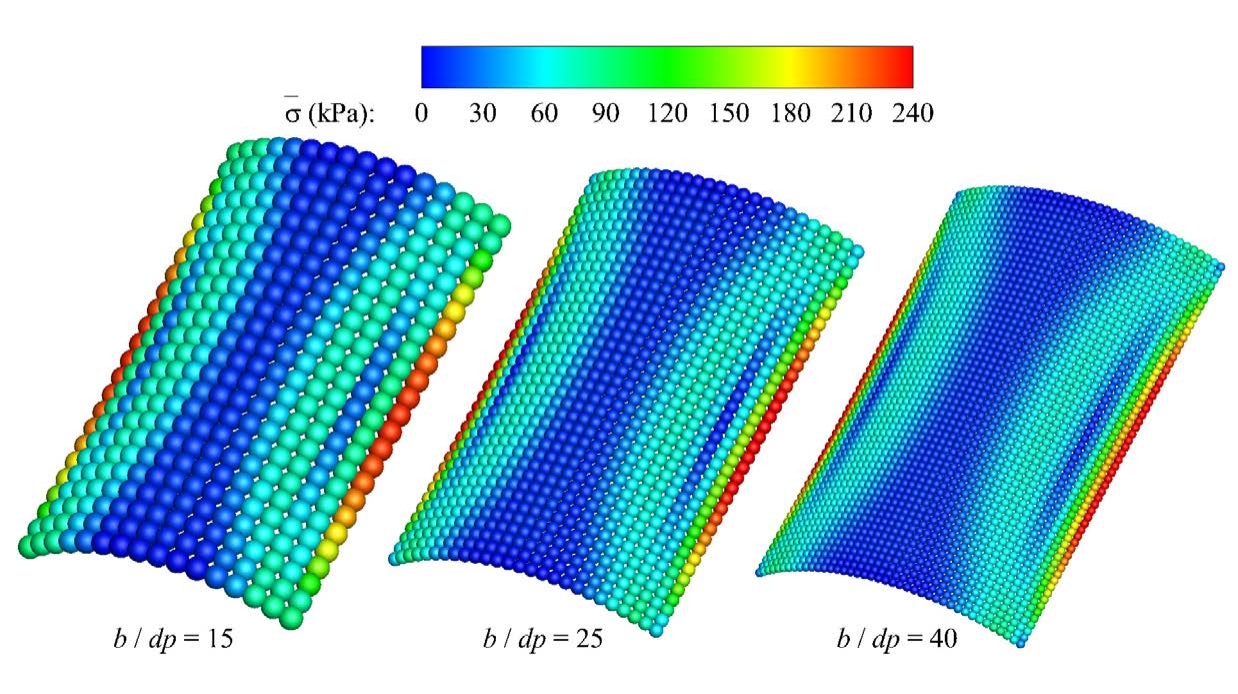}
	\caption{Scordelis-Lo roof: Particles colored by 
		the von Mises stress $\bar\sigma$ of the mid-surface 
		obtained by the present method with particle refinement. 
		The material is set as the density $\rho_0 = 36 ~\text{kg} / \text{m}^3$, 
		Young’s modulus $E = 432 ~\text{MPa}$ and Poisson’s ratio $\nu = 0.0$.}
	\label{figs:3D_roof_convergence_stress}
\end{figure}
\begin{figure}[htb!]
	\centering
	\includegraphics[trim = 2mm 6mm 2mm 2mm, width=0.5 \textwidth] {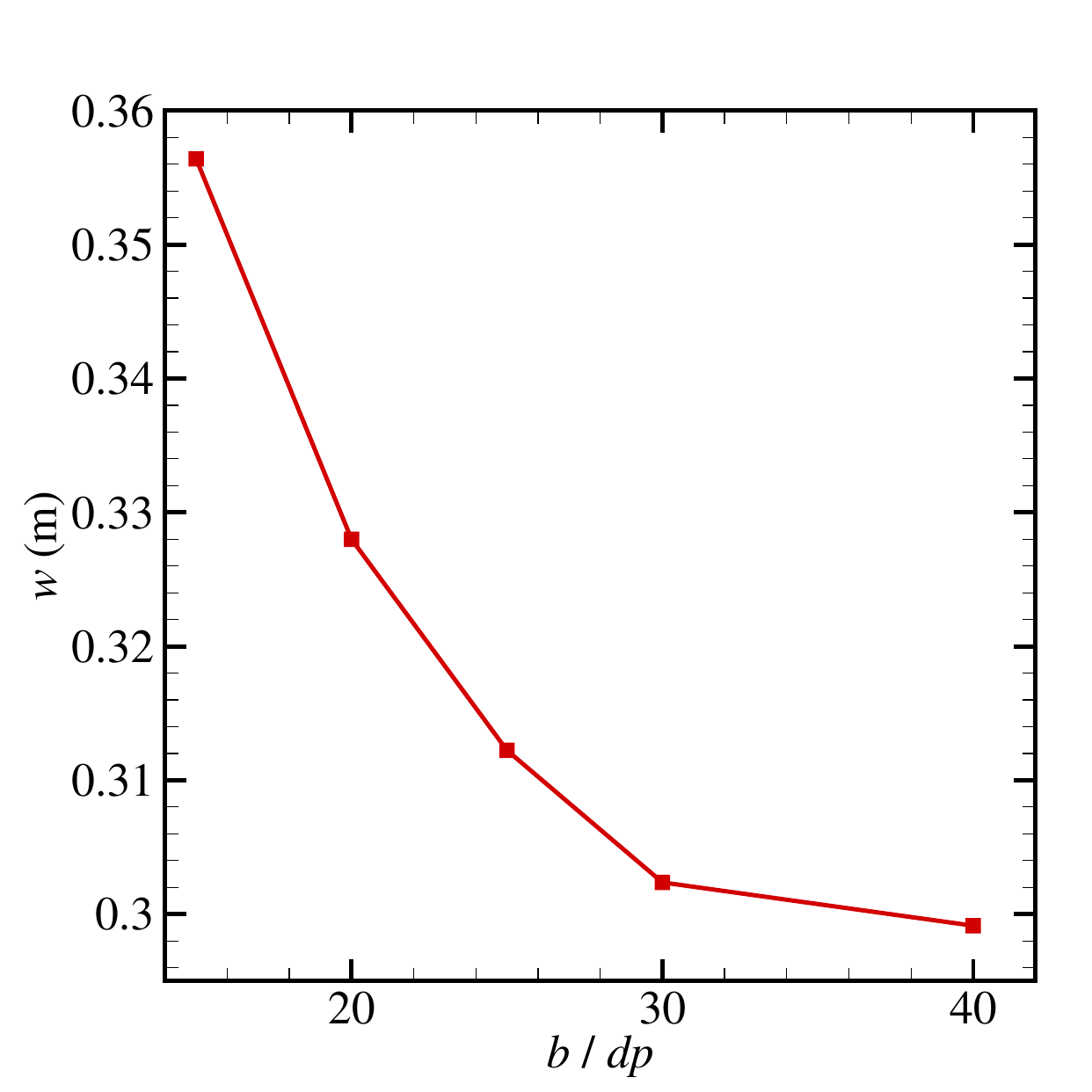}
	\caption{Scordelis-Lo roof: Convergence study of the displacement obtained by 
		the present method with particle refinement.}
	\label{figs:3D_roof_convergence}
\end{figure}
Figure \ref{figs:3D_roof_convergence_stress} shows 
the particles colored with the von Mises stress $\bar \sigma $ 
of the mid-surface
obtained at different resolutions. 
The regular particle distribution and smooth stress fields are noted.
With increasing resolution, 
a clear convergence is exhibited. 
The profile of displacement $w$ with varying spatial resolution 
obtained by the present method is depicted 
in Figure \ref{figs:3D_roof_convergence}.
It can be noted that the result converges rapidly to $w = 0.2991~\text{m}$ 
when $b/dp = 40$ with 1.09\% error compared 
to the solution of Refs. \cite{belytschko1985stress, simo1989stress}. 
\subsection{Pinched hemispherical shell}
\begin{figure}
	\begin{center}
		\includegraphics[trim = 2mm 6mm 2mm 2mm, width=\textwidth]{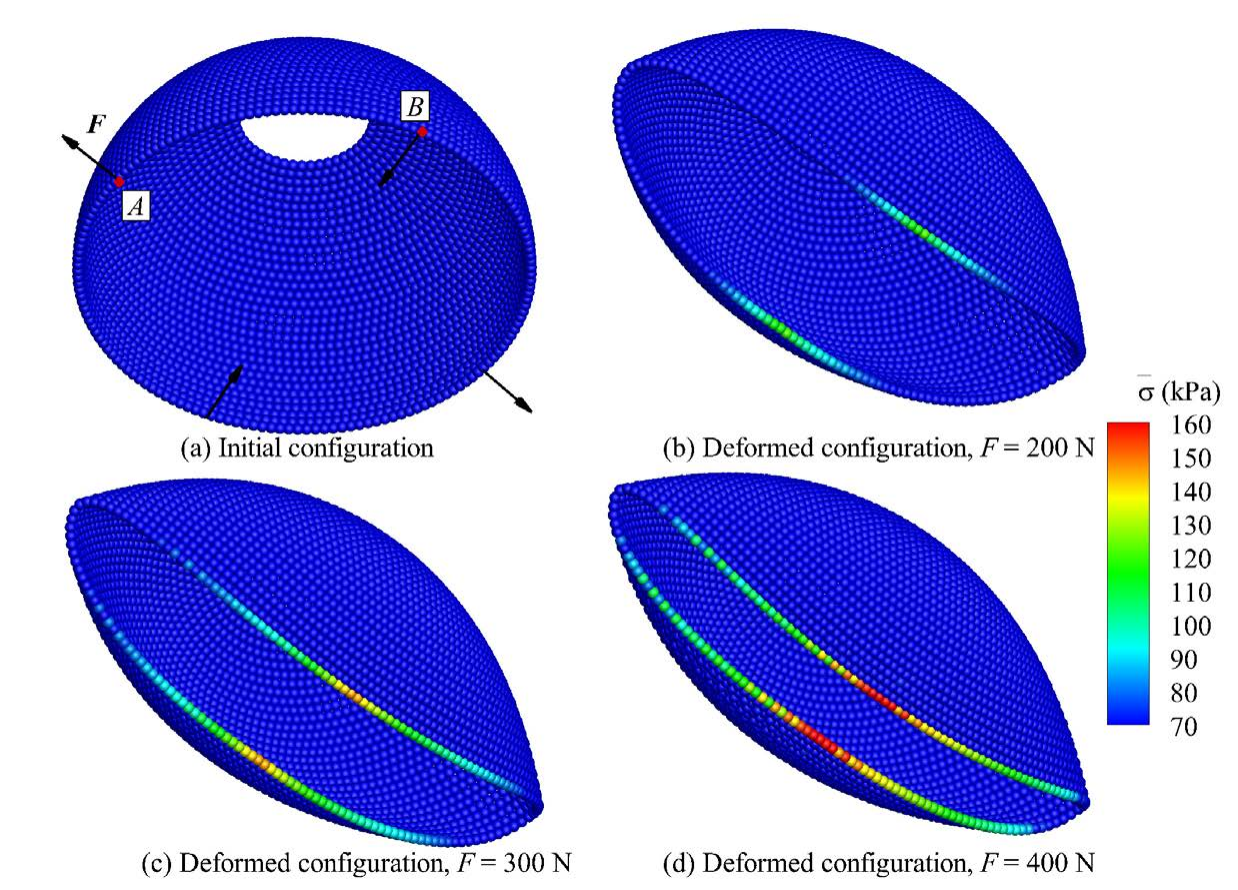}
		\caption{\small{Pinched hemispherical shell: 
				(a) Initial configuration with the radius of the mid-surface 
				$r = 10.0~\text{m}$ and thickness $d = 0.04~\text{m}$, 
				(b-d) particles colored by the von Mises stress $\bar\sigma$ of the mid-surface under 3 point force magnitudes
				at spatial discretication $2 \pi r / dp = 160$. 
				The material parameters are set as 
				the density $\rho_0 = 1100 ~\text{kg} / \text{m}^3$,  
				Young’s modulus $E = 68.25 ~\text{MPa}$ 
				and Poisson’s ratio $\nu = 0.3$.}} 
		\label{figs:3D_pinched_hemisphere}
	\end{center}
\end{figure}

We now consider a pinched hemispherical shell 
with an $18^\circ$ circular cutout at its pole
following Refs. \cite{simo1990stressIV, buechter1992shell, 
	jiang1994simple, sze2004popular, payette2014seven}. 
As shown in Figure \ref{figs:3D_pinched_hemisphere}(a), 
the hemispherical shell with the radius $r = 10.0~\text{m}$ 
and thickness $d = 0.04~\text{m}$ 
is loaded by four alternating radial point forces $\bm{F}$, 
prescribed along the equator at $90^\circ$ intervals. 
A linear elastic material with the density $\rho_0 = 1100 ~\text{kg} / \text{m}^3$,  
Young’s modulus $E = 68.25 ~\text{MPa}$ and Poisson’s ratio $\nu = 0.3$
is applied. 

\begin{figure}
	\begin{center}
		\includegraphics[trim = 2mm 6mm 2mm 2mm, width=\textwidth]{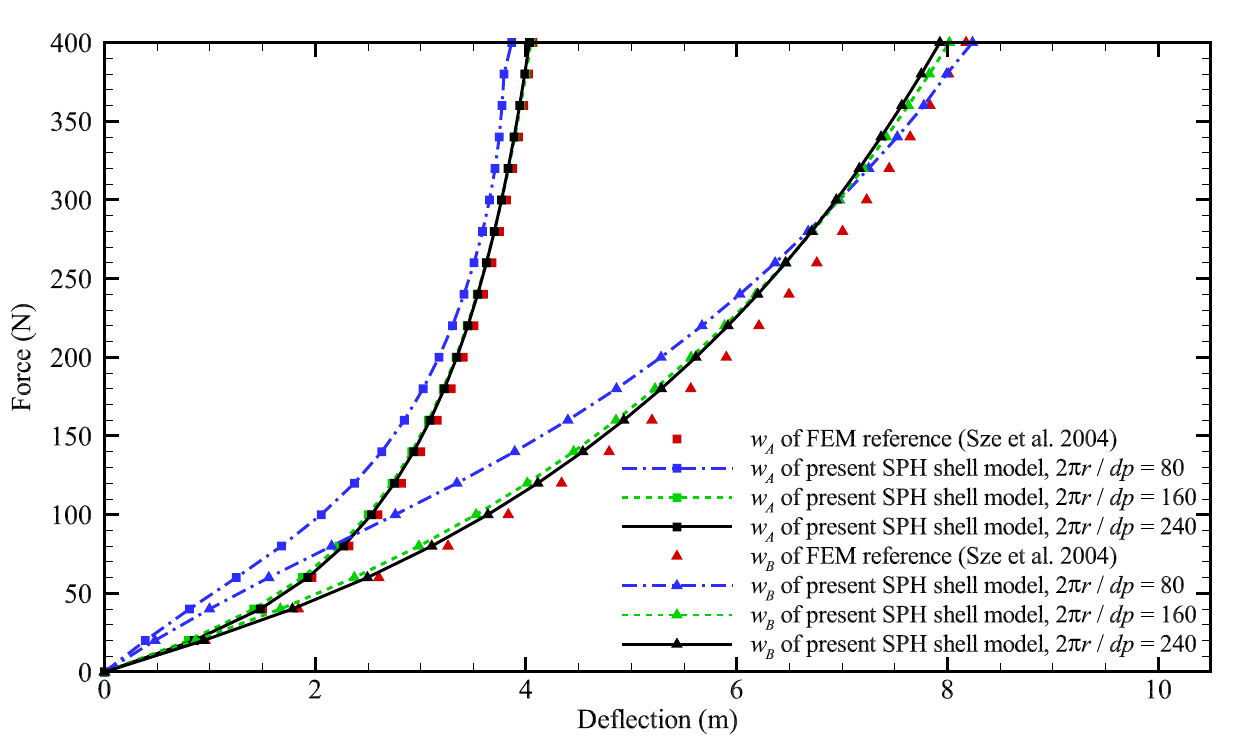}
		\caption{Pinched hemispherical shell: Curves of radical displacements of points $A$ and $B$ 
			as a function of the point force magnitude and spatial resolution, 
			and their comparison with those of  Sze et al. \cite{sze2004popular}.}
		\label{figs:3D_hemisphere_comparison}
	\end{center}
\end{figure}

Figure \ref{figs:3D_pinched_hemisphere}(b-d) shows 
the distribution of von Mises stress $\bar \sigma $ at the mid-surface 
under varying magnitude of the point force $\bm{F}$ . 
The regular particle distribution is observed, 
although slight stress fluctuation is exhibited near the place 
where the point force is applied. 
For quantitative analysis and convergence study, 
the radial deflections $w_A$ and $w_B$ of monitoring points $A$ and $B$ 
as a function of the point force magnitude and resolution 
are compared with those of Ref. \cite{sze2004popular}. 
Three different spatial discretizations, 
i.e., $2 \pi r/dp = 80, 160~\text{and}~240$,
are considered for convergence study. 
As shown in Figure \ref{figs:3D_hemisphere_comparison}, 
the results of present SPH shell model is quickly converging to 
those of Ref. \cite{sze2004popular}.
\subsection{Pulled-out cylindrical shell}
\begin{figure}
	\begin{center}
		\includegraphics[trim = 2mm 6mm 2mm 2mm, width=0.5\textwidth]{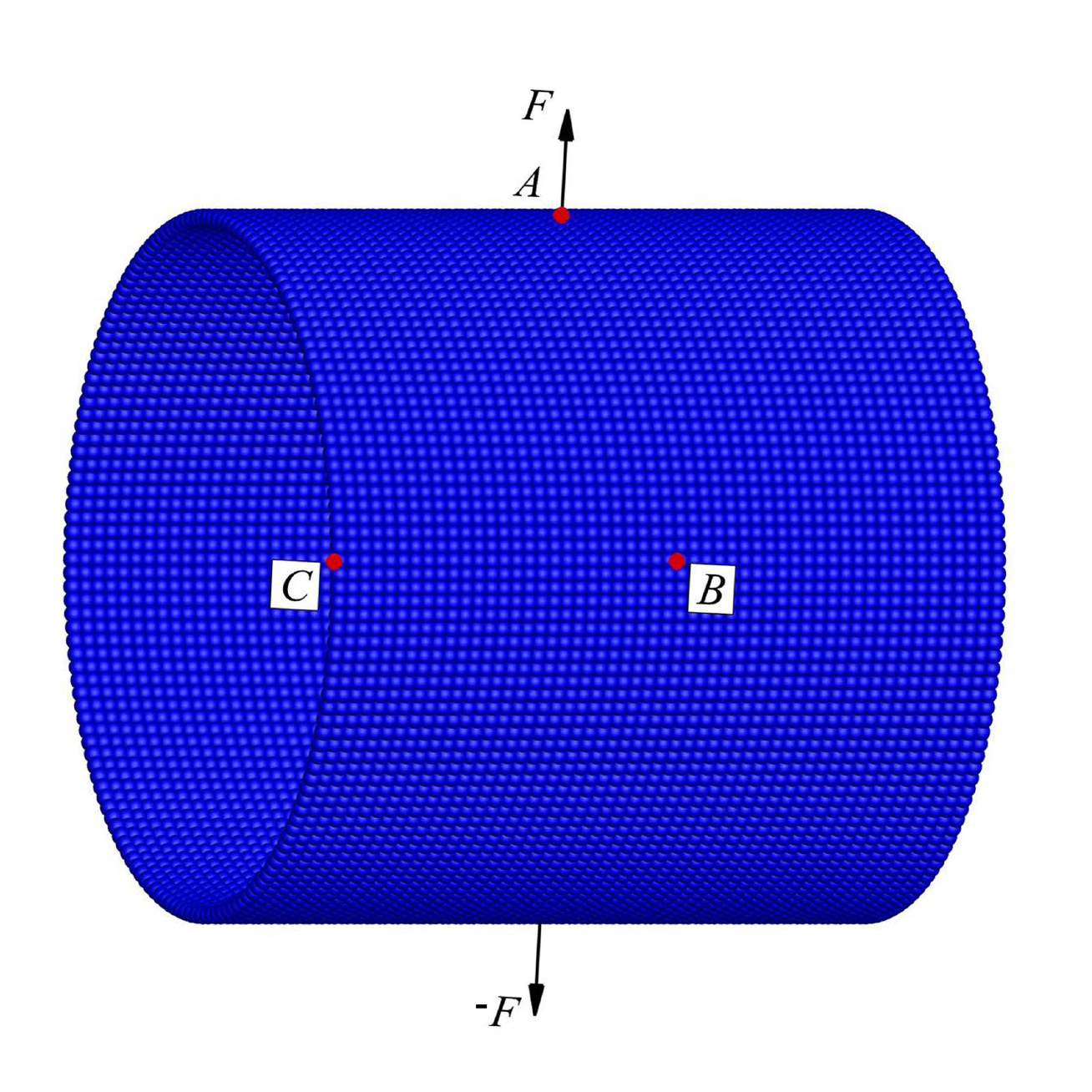}
		\caption{\small{Pulled-out cylindrical shell: Initial configuration with the radius of the mid-surface 
				$r = 5.0~\text{m}$, length $a = 10.35~\text{m}$ and thickness $d = 0.094~\text{m}$.}} 
		\label{figs:3D_cylinder_setup}
	\end{center}
\end{figure}

A more challenging benchmark test with large displacements 
is considered in this section following Refs. \cite{maurel2008sph, jiang1994corotational}.
As shown in Figure \ref{figs:3D_cylinder_setup}, 
a cylindrical shell with the radius $r = 5.0~\text{m}$, 
length $a = 10.35~\text{m}$ and thickness $d = 0.094~\text{m}$ 
is subjected to a pair of point forces $\bm{F}$ 
which are equal in magnitude and opposite in direction. 
A linear elastic material with the density $\rho_0 = 1100 ~\text{kg} / \text{m}^3$,  
Young’s modulus $E = 10.5 ~\text{MPa}$ and Poisson’s ratio $\nu = 0.3125$
is applied. 

\begin{figure}
	\begin{center}
		\includegraphics[trim = 2mm 6mm 2mm 2mm, width=\textwidth]{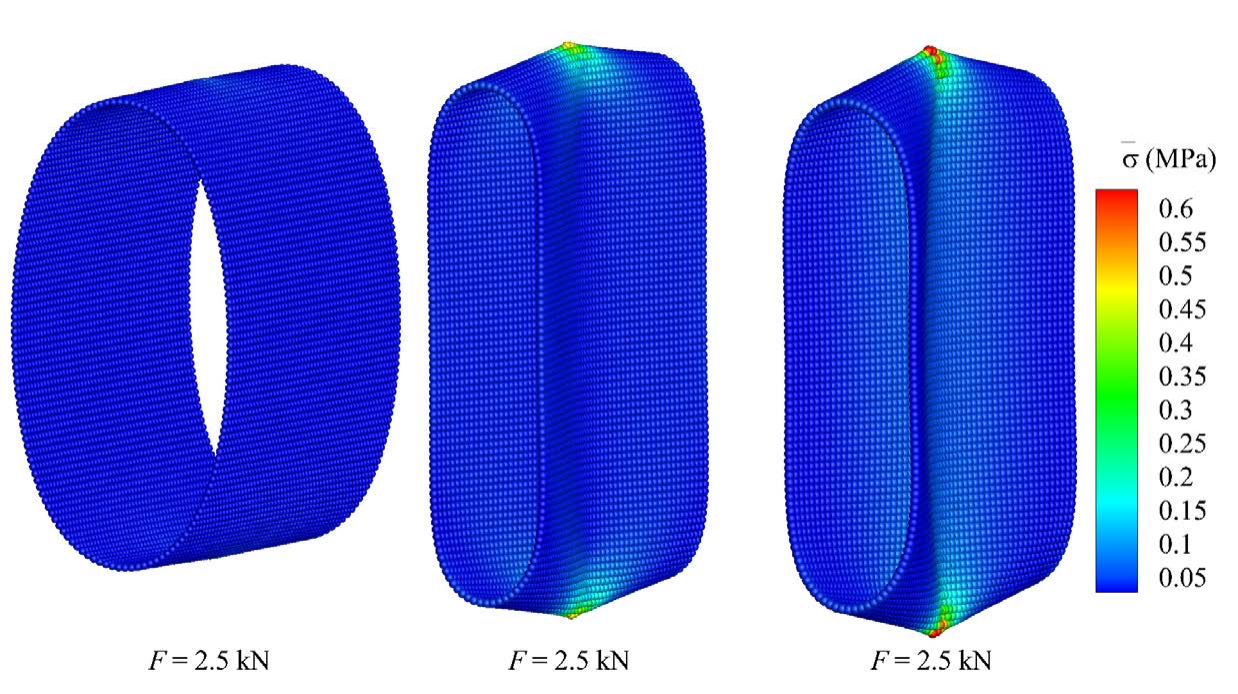}
		\caption{Pulled-out cylindrical shell: 
			Particles colored by the von Mises stress $\bar\sigma$ of the mid-surface
			under 3 point force magnitudes at spatial resolution of $b/dp = 240$. 
			The material parameters are set as the density $\rho_0 = 1100 ~\text{kg} / \text{m}^3$,  
			Young’s modulus $E = 10.5 ~\text{MPa}$ and Poisson’s ratio $\nu = 0.3125$.} 
		\label{figs:3D_cylinder_stress}
	\end{center}
\end{figure}

Figure \ref{figs:3D_cylinder_stress} shows 
the distribution of von Mises stress $\bar \sigma $ at the mid-surface 
under varying magnitude of the point force $\bm{F}$. 
The regular particle distribution and smooth stress fields, 
even close to the place where the point force is applied, 
are observed. 
\begin{figure}
	\begin{center}
		\includegraphics[trim = 2mm 6mm 2mm 2mm, width=\textwidth]{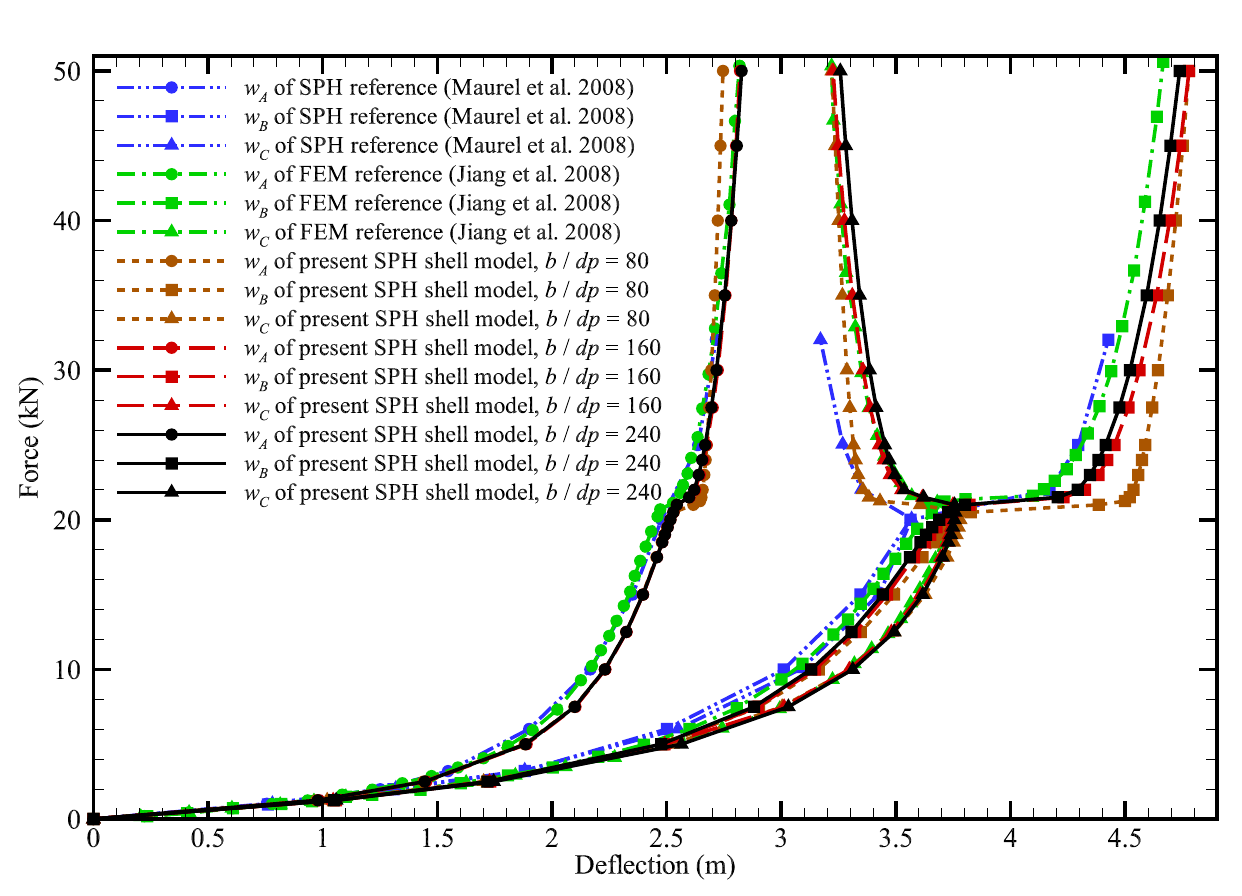}
		\caption{Pulled-out cylindrical shell: 
			Curves of radical displacements of points $A$, $B$ and $C$ 
			as a function of the point force magnitude and spatial resolution, 
			and their comparison with those of  Maurel and Combescure \cite{maurel2008sph} and Jiang et al. \cite{jiang1994corotational}.}
		\label{figs:3D_cylinder_comparison}
	\end{center}
\end{figure}
For quantitative analysis and convergence study, 
the radial displacements $w_A$, $w_B$ and $w_C$
of monitoring points $A$, $B$ and $C$ 
as a function of the point force magnitude and resolution 
are compared with those of Ref. \cite{maurel2008sph, jiang1994corotational}. 
Three different spatial discretizations, 
i.e., $b/dp = 80, 160~\text{and}~240$ with $b = 2 \pi r$ denoting 
the circumference length of the end, 
are considered for convergence study. 
As shown in Figure \ref{figs:3D_cylinder_comparison}, 
the bifurcation point of the curve is accurately predicted, 
suggesting good accuracy and robustness of the present method.
\subsection{Pinched semi-cylindrical shell}
\begin{figure}
	\begin{center}
		\includegraphics[trim = 2mm 6mm 2mm 2mm, width=\textwidth]{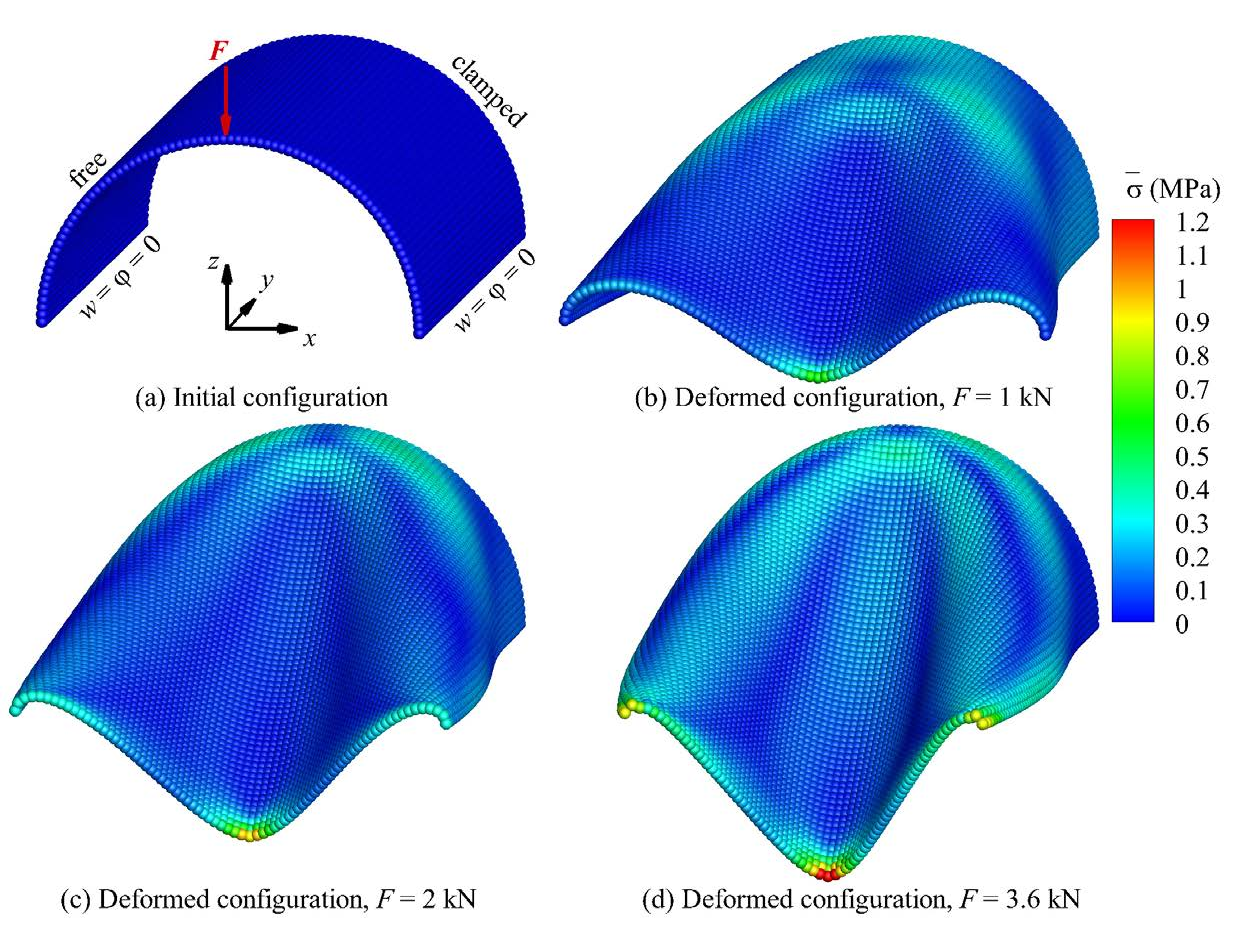}
		\caption{\small{Pinched semi-cylindrical shell: 
				(a) Initial configuration with the radius $r = 1.016~\text{m}$, 
				length $a = 3.048~\text{m}$ and thickness $d = 0.03~\text{m}$, 
				(b-d) particles colored by the von Mises stress $\bar\sigma$ of the mid-surface under 3 point force magnitudes 
				at spatial resolution $b / dp = 80$. 
				The material parameters are set as 
				the density $\rho_0 = 1100 ~\text{kg} / \text{m}^3$,  
				Young’s modulus $E = 20.685 ~\text{MPa}$ 
				and Poisson’s ratio $\nu = 0.3$.}} 
		\label{figs:3D_pinched_cylinder}
	\end{center}
\end{figure}

We further consider a pinched semi-cylindrical shell 
with finite deformation and rotation 
following Refs. \cite{stander1989assessment, brank1995implementation, 
	sze2004popular, arciniega2007tensor}. 
As shown in Figure \ref{figs:3D_pinched_cylinder}(a), 
the semi-cylindrical shell with the radius $r = 1.016~\text{m}$, 
length $a = 3.048~\text{m}$ and thickness $d = 0.03~\text{m}$ 
is completely clamped at a circumferential periphery 
and experiences a pinching force at the center of free-hanging periphery.
Along its longitudinal edges, 
the vertical direction and the rotation about the $y$-axis are constrained.
The elastic material properties are 
density $\rho_0 = 1100 ~\text{kg} / \text{m}^3$,  
Young’s modulus $E = 20.685 ~\text{MPa}$ and Poisson’s ratio $\nu = 0.3$. 

\begin{figure}
	\begin{center}
		\includegraphics[trim = 2mm 6mm 2mm 2mm, width=0.5\textwidth]{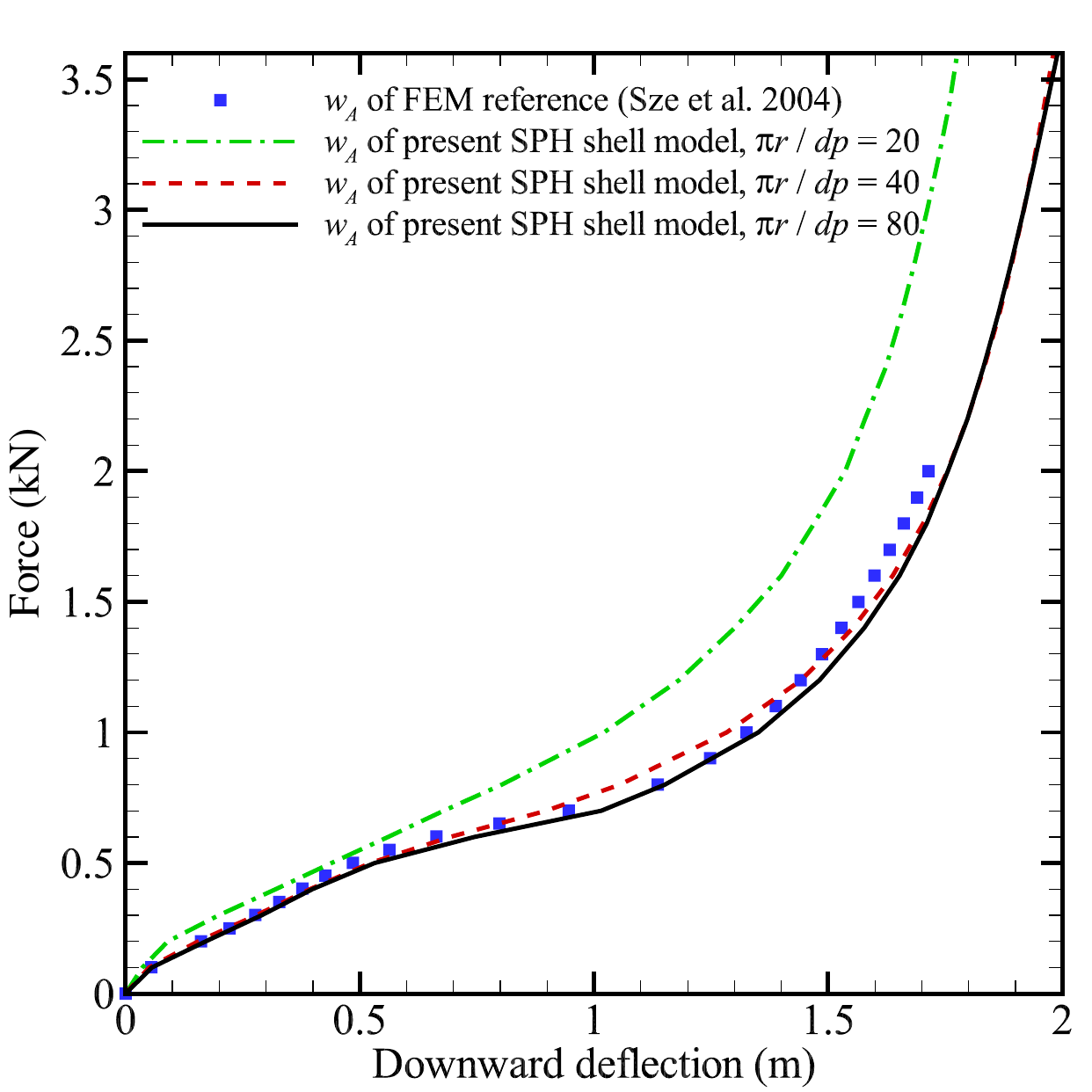}
		\caption{Pinched semi-cylindrical shell: 
			Curves of radical displacements of point $A$  
			as a function of the point force magnitude and spatial resolution, 
			and their comparison with those of  Sze et al. \cite{sze2004popular}.}
		\label{figs:3D_pinched_cylinder_comparison}
	\end{center}
\end{figure}

Figure \ref{figs:3D_pinched_cylinder}(b-d) shows 
the distribution of von Mises stress $\bar \sigma $ at the mid-surface 
under varying magnitude of the point force $\bm{F}$ . 
Noted that 
the present method features regular particle distribution and smooth stress fields, 
even close to the constrained edges and place where the point force is applied, 
without singularities for finite rotations (more than $0.5\pi$). 
For quantitative analysis and convergence study, 
the downward deflection $w_A$ of monitoring point $A$ 
as a function of the point force magnitude and resolution 
is compared with that of Ref. \cite{sze2004popular}. 
Three different spatial discretizations, 
i.e., $\pi r / dp = 20, 40~\text{and}~80$, 
are considered for convergence study. 
As shown in Figure \ref{figs:3D_pinched_cylinder_comparison}, 
the result difference obtained by the present SPH shell method 
between different resolution rapidly decreases
as the spatial refinement,  
and the results agree well with those of Ref. \cite{sze2004popular}.

\section{Concluding remarks}\label{sec:conclusion}
In this paper, 
we present a reduced-dimensional SPH method
for quasi-static and dynamic analyses 
of both thin and moderately thick plate and shell structures. 
By introducing two reduced-dimensional linear-reproducing correction matrices,
the method reproduces linear gradients of the position and pseudo-normal. 
The finite deformation is taken into account 
by considering all terms of strain 
with the help of Gauss-Legendre quadrature along the thickness. 
To cope with large rotations, 
the method introduces weighted non-singularity conversion relation 
between the rotation angles and pseudo normal.
A non-isotropic Kelvin-Voigt damping 
and a momentum-conserving hourglass control algorithm with a limiter are 
also proposed to increase numerical stability and to suppress hourglass modes. 
An extensive set of numerical examples have been investigated to
demonstrate the accuracy and robustness of the present method.
Note that, while the plate and shell structure considered 
here have moderate and high modulus,
one extension of the present method is for soft thin structures, such as membranes.    
Another outlook, 
along with the multi-physical modeling within unified computational framework,
is to develop the SPH method for the interaction between fluid and thin structures. 

\printnomenclature 

\clearpage
\section*{CRediT authorship contribution statement}
{\bfseries  Dong Wu:} Conceptualization, Methodology, Investigation, Visualization, Validation, Formal analysis, Writing - original draft, Writing - review \& editing; 
{\bfseries  Chi Zhang:} Investigation, Writing - review \& editing;
{\bfseries  Xiangyu Hu:} Supervision, Methodology, Investigation, Writing - review \& editing.

\section*{Declaration of competing interest }
The authors declare that they have no known competing financial interests 
or personal relationships that could have appeared to influence the work reported in this paper.

\section*{Acknowledgment}
D. Wu is partially supported by the China Scholarship Council (No. 201906130189). 
D. Wu, C. Zhang and X.Y. Hu would like to express their gratitude to the German Research Foundation (DFG) 
for its sponsorship of this research under grant number DFG HU1527/12-4.
\clearpage

\section*{References}
\vspace{-1.0cm}
\renewcommand{\refname}{}
\bibliographystyle{elsarticle-num}
{\small
	\bibliography{IEEEabrv,mybibfile}
}

%
%
\end{document}